\documentclass[12pt,oneside,reqno]{amsart}
\usepackage{mathrsfs}
\usepackage{graphics}
\usepackage{amssymb}
\usepackage{enumerate}
\newcommand{\dif}{\mathrm{d}}

\newcommand{\be}{\begin{eqnarray}}
\newcommand{\ee}{\end{eqnarray}}
\newcommand{\ce}{\begin{eqnarray*}}
\newcommand{\de}{\end{eqnarray*}}
\newtheorem{theorem}{Theorem}[section]
\newtheorem{lemma}[theorem]{Lemma}
\newtheorem{remark}[theorem]{Remark}
\newtheorem{definition}[theorem]{Definition}
\newtheorem{proposition}[theorem]{Proposition}
\newtheorem{Example}[theorem]{Example}
\newtheorem{corollary}[theorem]{Corollary}

\def\s{\sigma}

\def\a{\alpha}

\def\d{\delta}
\def\p{\partial}

\def\[{{\Big[}}
\def\]{{\Big]}}
\def\<{{\langle}}
\def\>{{\rangle}}
\def\({{\Big(}}
\def\){{\Big)}}

\def\no{\nonumber}
\def\bt{\begin{theorem}}
\def\et{\end{theorem}}
\def\bl{\begin{lemma}}
\def\el{\end{lemma}}
\def\br{\begin{remark}}
\def\er{\end{remark}}
\def\bx{\begin{Example}}
\def\ex{\end{Example}}
\def\bd{\begin{definition}}
\def\ed{\end{definition}}
\def\bp{\begin{proposition}}
\def\ep{\end{proposition}}
\def\bc{\begin{corollary}}
\def\ec{\end{corollary}}

\def\cD{{\mathcal D}}

\def\cK{{\mathcal K}}
\def\cL{{\mathcal L}}

\def\cO{{\mathcal O}}
\def\cP{{\mathcal P}}

\def\mE{{\mathbb E}}

\def\mN{{\mathbb N}}

\def\mP{{\mathbb P}}

\def\mR{{\mathbb R}}

\def\mW{{\mathbb W}}

\def\sA{{\mathscr A}}
\def\sB{{\mathscr B}}
\def\sC{{\mathscr C}}

\def\sF{{\mathscr F}}

\def\sL{{\mathscr L}}

\def\sV{{\mathscr V}}

\def\geq{\geqslant}
\def\leq{\leqslant}

\begin{document}

\allowdisplaybreaks

\title{Stability for multivalued McKean-Vlasov stochastic differential equations}

\author{Huijie Qiao$^*$  and  Jun Gong}



\thanks{This work was supported by NSF of China (No. 12071071).}

\thanks{$*$ Corresponding author: hjqiaogean@seu.edu.cn}


\date{}
\dedicatory{School of Mathematics,
Southeast University\\
Nanjing, Jiangsu 211189, China}

\begin{abstract}
The work concerns multivalued McKean-Vlasov stochastic differential equations. First of all, we prove the existence and uniqueness of strong solutions for multivalued McKean-Vlasov stochastic differential equations with non-Lipschitz coefficients. Then, the classical It\^{o}'s formula is extended to that for multivalued McKean-Vlasov stochastic differential equations. Finally, the asymptotic stability of second moments and the almost surely asymptotic stability for their solutions in terms of a Lyapunov function are shown. 
\end{abstract}

\maketitle \rm

\section{Introduction}

Given a filtered probability space $(\Omega,\mathscr{F},\{\mathscr{F}_t\}_{t\geq 0},\mP)$ and a $m$-dimensional standard Brownian motion $W_{\cdot}=(W_{\cdot}^1,W_{\cdot}^2,\cdots, W_{\cdot}^m)$ defined on it. Consider the following multivalued McKean-Vlasov stochastic differential equation (SDE for short) on $\mR^d$:
\be\left\{\begin{array}{l}
\dif X_t\in \ -A(X_t)\dif t+ \ b(X_t,\sL_{X_t})\dif t+\sigma(X_t,\sL_{X_t})\dif W_t,\\
X_0=\xi,     
\end{array}
\label{eq1}
\right.
\ee
where $\xi$ is a $\sF_0$-measurable random variable with $\mE|\xi|^{2}<\infty$, $A:\mR^d \mapsto 2^{\mR^d} $ is a maximal monotone operator, $\sL_{X_t}$ is the probability  distribution of $X_t$ with respect to the probability measure $\mP$, and the coefficients $b:\mR^d\times\cP_2(\mR^d)\mapsto{\mR^d}, \,\,\sigma:\mR^d\times\cP_2(\mR^d)\mapsto{\mR^d}\times{\mR^m}$ are Borel measurable ($\cP_2(\mR^d)$ is defined in Subsection \ref{nota}).

If $A = 0$, Eq.(\ref{eq1}) becomes McKean-Vlasov SDE. And the first work on McKean-Vlasov SDEs can be tracked back to McKean \cite{MC}, who was inspired of Kac's programme in Kinetic theory in \cite{ka}. From then on, a large number of results appear (c.f. \cite{LM}-\cite{WDL}). Let us mention some related works. For Eq.(\ref{eq1}), under non-Lipschitz coefficients, Ding and Qiao \cite{DQ1, DQ2} proved the well-posedness, exponential stability of the second moment, exponentially 2-ultimate boundedness and almost surely asymptotic stability of strong solutions to these equations. In \cite{WDL}, Hammersley, Siska and Szpruch also showed exponentially 2-ultimate boundedness of strong solutions for Eq.(\ref{eq1}) under the local boundedness for $b, \sigma$ and some extra properties for integrated Lyapunov functions. 

If $A \ne 0$, Eq.(\ref{eq1}) is called a multivalued McKean-Vlasov SDE. Under global Lipschitz conditions, Chi \cite{CHI} proved the existence and uniqueness of strong solutions for the following multivalued McKean-Vlasov SDE on $\mR^d$:
\be
\dif X_t\in -A(X_t)\dif t+\bar{b}[X_t,\sL_{X_t}]\dif t+\bar{\sigma}[X_t,\sL_{X_t}]\dif W_t,
\label{mv2}
\ee
where for any $x\in\mR^d$ and $\mu\in\cP_2(\mR^d)$
$$
\bar{b}[x,\mu]:=\int_{\mR^d}\bar{b}(x,y)\mu(\dif y), \quad \bar{\sigma}[x,\mu]:=\int_{\mR^d}\bar{\sigma}(x,y)\mu(\dif y),
$$
and $\bar{b}: \mR^d\times\mR^d\mapsto\mR^d$, $\bar{\sigma}: \mR^d\times\mR^d\mapsto\mR^d\times\mR^m$ are Borel measurable. When the operator $A$ is the sub-differential of some convex function, Ren and Wang \cite{RW}  studied the following multivalued McKean-Vlasov SDE on $\mR^d$:
\be
 \dif X_t\in -A(X_t)\dif t+\tilde{b}(X_t,\mE[X_t])\dif t+\tilde{\sigma}(X_t,\mE[X_t])\dif W_t,
\label{mv3}
\ee
where $\tilde{b}: \mR^d\times\mR^d\mapsto\mR^d$, $\tilde{\sigma}:\mR^d\times\mR^d\mapsto\mR^d\times\mR^m$ are Borel measurable and $\mE[X_t]$ denotes the expectation with respect to the probability measure $\mP$. If $\tilde{b}, \tilde{\s}$ satisfy global Lipschitz conditions, they showed well-posedness and a large deviation principle for Eq.(\ref{mv3}). Note that Eq.(\ref{eq1}) is more general than Eq.(\ref{mv2}) and (\ref{mv3}). However, as far as we know, Eq.(\ref{eq1}) seems not to be studied in the literature.

In this paper, we first study the existence and uniqueness of strong solutions to Eq.(\ref{eq1}) under non-Lipschitz conditions. Our result (Theorem \ref{EU}) can cover \cite[Theorem 3.1]{CHI} and \cite[Theorem 3.3]{RW} for the time homogeneous case. Then, the classical It\^{o} formula for SDEs is extended to that for multivalued McKean-Vlasov SDEs. We emphasize that the It\^{o} formula for multivalued McKean-Vlasov SDEs (Theorem \ref{if}) is different from the It\^{o} formula for McKean-Vlasov SDEs (c.f. \cite[Proposition 3.1]{QW}). After this, the stability for the strong solution to Eq.$(\ref{eq1})$ is considered. We offer sufficient conditions to assure the asymptotic stability of the second moment in terms of Lyapunov functions. Finally, we prove the almost surely asymptotic stability for the strong solution to Eq.$(\ref{eq1})$. Here we remind that the appearance of the maximal monotone operator causes a lot of trouble in the deduction.

The rest of the paper is organized as follows. In Section \ref{fram}, we recall some basic notation and introduce maximal monotone operators and derivatives for functions on $\cP_{2}(\mR^d)$ with respect to the measures. In Section \ref{SU}, we prove that Eq.(\ref{eq1}) has a unique strong solution under non-Lipschitz conditions. Next we extend the classical It\^o formula to that for multivalued McKean-Vlasov SDEs in Section \ref{geneito}. Then in Section \ref{stab} we present the asymptotic stability of second moments and the almost surely asymptotic stability for the strong solution to Eq.$(\ref{eq1})$. Finally, an example is given to explain our result in Section \ref{app}.

The following convention will be used throughout the paper: $C$ with or without indices will denote different positive constants whose values may change from one place to another.

\section{Preliminary}\label{fram}

In the section, we introduce notations and concepts and recall some results used in the sequel.

\subsection{Notations}\label{nota}

In the subsection, we introduce some notations.

For convenience, we shall use $\mid\cdot\mid$ and $\parallel\cdot\parallel$  for norms of vectors and matrices, respectively. Furthermore, let $\langle\cdot$ , $\cdot\rangle$ denote the scalar product in $\mR^d$. Let $B^*$ denote the transpose of the matrix $B$.

Let $C(\mR^d)$ be the collection of continuous functions on $\mR^d$ and $C^2(\mR^d)$ be the space of continuous functions on $\mR^d$ which have continuous partial derivatives of order up to $2$. 

Let $\sB(\mR^d)$ be the Borel $\sigma$-algebra on $\mR^d$ and $\cP({\mR^d})$ be the space of all probability measures defined on $\sB(\mR^d)$ carrying the usual topology of weak convergence. Let $\cP_{2}(\mR^d)$ be the set of probability measures on $\sB(\mR^d)$ with finite second order moments. That is,
$$
\cP_2\left( \mathbb{R}^d \right) :=\left\{ \mu \in \cP\left( \mathbb{R}^d \right): \mu(|\cdot|^{2}):=\int_{\mathbb{R}^d}{\left| x \right|^2\mu \left( \dif x \right) <\infty} \right\}. 
$$
It is known that $\cP_2(\mR^d)$ is a Polish space endowed with the $L^2$-Wasserstein distance defined by
$$
\mathbb{W}_2(\mu,\nu):= \inf\limits_{\pi\in\Psi(\mu,\nu)}\left(\int_{\mathbb{R}^d\times\mathbb{R}^d}|x-y|^{2}\pi(\dif x,\dif y)\right)^{\frac{1}{2}}, \quad \mu , \nu\in \cP_2(\mR^d),
$$
where $\Psi(\mu,\nu)$ is the set of all couplings $\pi$ with marginal distributions $\mu$ and $\nu$. Moreover, if $\xi,\zeta$ are two random variables with distributions $\sL_\xi, \sL_\zeta$ under $\mP$, respectively,
$$
\mathbb{W}_2(\sL_\xi, \sL_\zeta)\leq (\mE|\xi-\zeta|^2)^{\frac{1}{2}},
$$
where $\mE$ stands for the expectation with respect to $\mP$.

\subsection{Maximal monotone operators}

In the subsection, we introduce maximal monotone operators. 

Fix a multivalued operator $A: \mR^d\mapsto 2^{\mR^d}$, where $2^{\mR^d}$ stands for all the subsets of $\mR^d$, and set
$$
\cD(A):= \left\{x\in \mR^d: A(x) \ne \emptyset\right\}
$$
and
$$
Gr(A):= \left\{(x,y)\in \mR^{2d}:x \in \cD(A), ~ y\in A(x)\right\}.
$$
Then we say that $A$ is monotone if $\langle x_1 - x_2, y_1 - y_2 \rangle \geq 0$ for any $(x_1,y_1), (x_2,y_2) \in Gr(A)$, and $A$ is maximal monotone if 
$$
(x_1,y_1) \in Gr(A) \iff \langle x_1-x_2, y_1 -y_2 \rangle \geq 0, \quad \forall (x_2,y_2) \in Gr(A).
$$

\bx\label{consetex}
Let $\cO\subset\mR^d$ be a closed convex set with non-empty interior. The indicator function of $\cO$ is defined by
\ce
I_{\cO}(x):=\left\{\begin{array}{ll}0, \qquad x\in \cO,\\
+\infty, \quad x\notin \cO.
\end{array}
\right.
\de
Then one can justify that $I_{\cO}$ is a lower semicontinuous convex function. The sub-differential of $I_{\cO}$ is given by
\ce
\partial I_{\cO}(x):=\left\{\begin{array}{ll}
\emptyset, \qquad x\notin \cO,\\
\{0\}, \quad x\in {\rm Int}(\cO),\\
\Pi_x, \qquad x\in \partial\cO,
\end{array}
\right.
\de
where $\Pi_x$ is the cone of unit outward normal to $\cO$ at $x$. It is well known that $\partial I_{\cO}$ is a maximal monotone operator with ${\rm Int}(\cD(\partial I_{\cO}))={\rm Int}(\cO)\neq \emptyset$.
\ex

Given $T>0$. Let $\sV_{0}$ be the set of all continuous functions $K: [0,T]\mapsto\mR^{d}$ with finite variations and $K_{0} = 0$. For $K\in\sV_0$ and $s\in [0,T]$, we shall use $|K|_{0}^{s}$ to denote the variation of $K$ on $[0,s]$. Set
\ce
&&\sA:=\Big\{(X,K): X\in C([0,T],\overline{\cD(A)}), K \in \sV_0, \\
&&\qquad\qquad\quad~\mbox{and}~\langle X_{t}-x, \dif K_{t}-y\dif t\rangle \geq 0 ~\mbox{for any}~ (x,y)\in Gr(A)\Big\}.
\de
Then about $\sA$ we recall three following results. (c.f. \cite{cepaa,ZXCH})

\bl\label{equi}
For $X\in C([0,T],\overline{\cD(A)})$ and $K\in \sV_{0}$, the following statements are equivalent:
\begin{enumerate}[(i)]
\item $(X,K)\in \sA$.
\item For any $(x,y)\in C([0,T],\mR^d)$ with $(x_t,y_t)\in Gr(A)$, it holds that 
$$
\left\langle X_t-x_t, \dif K_t-y_t\dif t\right\rangle \geq0.
$$
\item For any $(X^{'},K^{'})\in \sA$, it holds that 
$$
\left\langle X_t-X_t^{'},\dif K_t-\dif K_t^{'}\right\rangle \geq0.
$$
\end{enumerate}
\el

\bl\label{inteineq}
Assume that ${\rm Int}(\cD(A))\ne\emptyset$, where ${\rm Int}(\cD(A))$ denotes the interior of the set $\cD(A)$. For any $a\in {\rm Int}(\cD(A))$, there exists constants $\gamma_1 >0$, and $\gamma_{2},\gamma_{3}\geq0$ such that  for any $(X,K)\in \sA$ and $0\leq s<t\leq T$,
$$
\int_s^t{\left< X_r-a, \dif K_r \right>}\geq \gamma_1\left| K \right|_{s}^{t}-\gamma _2\int_s^t{\left| X_r-a\right|}\dif r-\gamma_3\left( t-s \right) .
$$
\el

\bl\label{limiconv}
Assume that $\{K^{n},n\in \mN\}\subset \sV_{0}$ converges to some $K$ in $C([0,T];\mR^d)$ and $\underset{n\in \mathbb{N}}{\sup}\left| K^n \right|_0^T<\infty $. Then $K\in \sV_0$, and 
$$
\underset{n\rightarrow \infty}{\lim}\int_0^T{\left< X_{s}^{n}, \dif K_{s}^{n} \right>}=\int_0^T{\left< X_{s}, \dif K_{s}\right>},
$$
where the sequence $\{X^{n}\}\subset C([0,T];\mR^d)$ converges to some $X$ in $C([0,T];\mR^d)$.
\el

\subsection{Derivatives for functions on $\cP_{2}(\mR^d)$}

In the subsection, we introduce derivatives for functions on $\cP_{2}(\mR^d)$ (c.f. \cite{Lion}). 

A function $f:\cP_{2}(\mR^d)\mapsto\mR$ is differential at $\mu\in \cP_{2}(\mR^d)$, if for $\tilde{f}(\gamma):=f(\sL_\gamma),\gamma\in L^2(\Omega,\mathscr{F},\mP;\mR^d)$, there exists some $\zeta\in L^2(\Omega,\mathscr{F},\mP;\mR^d)$ with $\sL_\zeta=\mu$ such that $\tilde{f}$ is Fr\'echet differentiable at $\zeta$, that is, there exists a linear continuous mapping $D\tilde{f}(\zeta):L^2(\Omega,\mathscr{F},\mP;\mR^d)\mapsto\mR$ such that for any $\eta\in L^2(\Omega,\mathscr{F},\mP;\mR^d)$
\ce
\tilde{f}(\zeta+\eta)-\tilde{f}(\zeta)=D\tilde{f}(\zeta)(\eta)+o(|\eta|_{L^2}),  \quad|\eta|_{L^2}\rightarrow0.
\de
Since $D\tilde{f}(\zeta)\in L(L^2(\Omega,\mathscr{F},\mP;\mR^d),\mR)$, it follows from the Riesz representation theorem that there exists a $\mP$-a.s. unique variable $\vartheta\in L^2(\Omega,\mathscr{F},\mP;\mR^d)$ such that for all $\eta\in L^2(\Omega,\mathscr{F},\mP;\mR^d)$
\ce
D\tilde{f}(\zeta)(\eta)=(\vartheta,\eta)_{L^2}=\mE[\vartheta\cdot\eta].
\de

\bd\label{f1}
We say that $f\in C^1(\cP_{2}(\mR^d))$, if there exists for all $\gamma\in L^2(\Omega,\mathscr{F},\mP;\mR^d)$ a $\sL_{\gamma}$-modification of $\partial_\mu f(\sL_\gamma)(\cdot)$, again denoted by $\partial_\mu f(\sL_\gamma)(\cdot)$, such that $\partial_\mu f:\cP_{2}(\mR^d)\times\mR^d\mapsto\mR^d$ is continuous, and we identify this continuous function $\partial_\mu f$ as the derivative of $f$.
\ed

\bd\label{f2}
We say that $f\in C^2(\cP_{2}(\mR^d))$, if for any $\mu\in \cP_{2}(\mR^d)$, $f\in C^1(\cP_{2}(\mR^d))$ and $\partial_\mu f(\sL_\gamma)(\cdot)$ is differentiable, and its derivative $\partial_y\partial_\mu f:\cP_{2}(\mR^d)\times\mR^d\mapsto\mR^d\otimes\mR^d$ is continuous, and for any $y\in \mR^d$, $\partial_{\mu}f(\cdot)(y)$ is differentiable, and its derivative $\partial_{\mu}^{2}f:\cP_{2}(\mR^d)\times\mR^d\times\mR^d\mapsto\mR^d\otimes\mR^d$ is continuous.
\ed

\bd
A function $F: \mR^d\times\cP_{2}(\mR^d)\mapsto\mR$  is said to be in $C^{2,2}(\mR^d\times\cP_{2}(\mR^d))$, if 
$(i)$ $F$ is $C^2$ in $x\in\mR^d$ and $\mu\in\cP_2(\mR^d)$, respectively;\\
$(ii)$ these derivatives 
\ce
\partial_{x}F(x,\mu), \partial_{x}^2F(x,\mu), \partial_{\mu}F(x,\mu)(y), \partial_{y}\partial_{\mu}F(x,\mu)(y), \partial_{\mu}^2F(x,\mu)(y,y^{'})
\de
are jointly continuous in the variable family $(x,\mu), (x,\mu,y)$ and $(x,\mu,y,y^{'})$, respectively.
\ed

\bd
A function $F: \mR^d\times\cP_{2}(\mR^d)\mapsto\mR$ is said to be in $\sC(\mR^d\times\cP_2(\mR^d))$, if $F\in C^{2,2}(\mR^d\times\cP_2(\mR^d))$ and for any compact set $\cK\subset\mR^d\times\cP_2(\mR^d)$,
$$
\sup\limits_{(x,\mu)\in\cK}\int_{\mR^d}\left(\|\partial_y\partial_\mu F(x,\mu)(y)\|^2+|\partial_\mu F(x,\mu)(y)|^2\right)\mu(\dif y)<\infty.
$$
\ed

\bd
A function $F: \mR^d\times\cP_{2}(\mR^d)\mapsto\mR$ is said to be in $C_{b}^{2,2}(\mR^d\times\cP_{2}(\mR^d))$, if 
$(i)$ $F\in C^{2,2}(\mR^d\times\cP_{2}(\mR^d))$;\\
$(ii)$ $F(x,\mu)$ and all its derivatives
\ce
\partial_{x}F(x,\mu), \partial_{x}^2F(x,\mu), \partial_{\mu}F(x,\mu)(y), \partial_{y}\partial_{\mu}F(x,\mu)(y), \partial_{\mu}^2F(x,\mu)(y,y^{'})
\de
are uniformly bounded in the variable family $(x,\mu), (x,\mu,y)$ and $(x,\mu,y,y^{'})$, respectively.
\ed

\bd\label{r2}
The function $F: \mR^d\times\cP_{2}(\mR^d)\mapsto\mR$ is said to be in $C_b^{2,2;1}(\mR^d\times\cP_2(\mR^d))$, if 
$(i)$ $F\in C_{b}^{2,2}(\mR^d\times\cP_{2}(\mR^d))$;\\
$(ii)$ $F(x,\mu)$ and all its derivatives
\ce
\partial_{x}F(x,\mu), \partial_{x}^2F(x,\mu), \partial_{\mu}F(x,\mu)(y), \partial_{y}\partial_{\mu}F(x,\mu)(y), \partial_{\mu}^2F(x,\mu)(y,y^{'})
\de
are Lipschitz continuous in the variable family $(x,\mu), (x,\mu,y)$ and $(x,\mu,y,y^{'})$, respectively.\\
If $F\in C_b^{2,2;1}(\mR^d\times\cP_2(\mR^d))$ and $F\geq 0$, we say that $F\in C_{b,+}^{2,2;1}(\mR^d\times\cP_2(\mR^d))$.
\ed

\section{The existence and uniqueness of strong solutions}\label{SU}
In this section, we study the existence and uniqueness of strong solution for Eq.(\ref{eq1}). 

First of all, we define strong solutions, weak solutions and the pathwise uniqueness of weak solutions for Eq.(\ref{eq1}). Fix $T>0$ and consider Eq.(\ref{eq1}), i.e.
\ce\left\{\begin{array}{l}
\dif X_t\in -A(X_t)\dif t+b(X_t,\sL_{X_t})\dif t+\sigma(X_t,\sL_{X_t})\dif W_t, \quad 0\leq t\leq T,\\
X_0=\xi.
\end{array}
\right.
\de

\bd(Strong solutions)\label{strosolu}
We say that Eq.$(\ref{eq1})$ admits a strong solution with the initial value $\xi$ if there exists a pair of adapted processes $(X,K)$ on a filtered probability space $(\Omega, \mathscr{F}, \{\mathscr{F}_t\}_{t\in[0,T]}, \mP)$ such that

$(i)$ $\mP(X_0=\xi)=1$,

$(ii)$ $X_t\in{\mathscr{F}_t^W}$, where $\{\mathscr{F}_t^W\}_{t\in[0,T]}$ stands for the $\sigma$-field filtration generated by $W$,

$(iii)$ $(X_{\cdot}(\omega),K_{\cdot}(\omega))\in \sA$ a.s. $\mP$,

$(iv)$ it holds that
\ce
\mP\left\{\int_0^T(\mid{b(X_s,\sL_{X_s})}\mid+\parallel{\sigma(X_s,\sL_{X_s})}\parallel^2)\dif s<+\infty\right\}=1,
\de
and
\ce
X_t=\xi-K_{t}+\int_0^tb(X_s,\sL_{X_s})\dif s+\int_0^t\sigma(X_s,\sL_{X_s})\dif W_s, \quad 0\leq{t}\leq{T}. 
\de
\ed

From the above definition, we know that $\sL_{X_0}=\sL_{\xi}$.

\bd(Weak solutions)\label{weaksolu}
We say that  Eq.(\ref{eq1}) admits a weak solution with the initial law $\sL_{\xi}\in \cP(\mR^{d})$, if there exists a filtered probability space $(\hat{\Omega}, \hat{\sF}, \{\hat{\sF}_t\}_{t\in[0,T]}, \hat{\mP})$, a $d$-dimensional standard $(\hat{\sF}_{t})$-Brownian motion $\hat{W}$ as well as a pair of $(\hat{\sF}_{t})$-adapted process $(\hat{X},\hat{K})$ defined on $(\hat{\Omega}, \hat{\sF}, \{\hat{\sF}_t\}_{t\in[0,T]}, \hat{\mP})$ such that
\begin{enumerate}[(i)]
	\item $\hat{\mP}\circ\hat{X}^{-1}_0=\sL_{\xi}$,
	\item $(\hat{X}_{\cdot}(\omega),\hat{K}_{\cdot}(\omega))\in \sA$ a.s. $\hat{\mP}$,
	\item it holds that  
	$$
\mP\left\{\int_0^T(\mid{b(\hat{X}_{s},\sL_{\hat{X}_{s}})}\mid+\parallel{\sigma(\hat{X}_{s},\sL_{\hat{X}_{s}})}\parallel^2)\dif s<+\infty\right\}=1,
	$$
	and
	$$
	\hat{X}_{t}=\hat{X}_{0}-\hat{K}_{t}+\int_{0}^{t}b(\hat{X}_{s},\sL_{\hat{X}_{s}})\dif s +\int_{0}^{t}\sigma(\hat{X}_{s},\sL_{\hat{X}_{s}})\dif \hat{W}_{s},\quad 0\leq{t}\leq{T}.
	$$
\end{enumerate}
\ed

Such a solution will be denoted  by $(\hat{\Omega}, \hat{\sF}, \{\hat{\sF}_t\}_{t\in[0,T]}, \hat{\mP}; \hat{W}, (\hat{X},\hat{K}))$.

\bd(Pathwise Uniqueness) \label{pathuniq}
Suppose $(\hat{\Omega}, \hat{\sF}, \{\hat{\sF}_t\}_{t\in[0,T]}, \hat{\mP}; \hat{W}, (\hat{X}^1,\hat{K}^1))$ and $(\hat{\Omega}, \hat{\sF}, \{\hat{\sF}_t\}_{t\in[0,T]}, \hat{\mP}; \hat{W}, (\hat{X}^2,\hat{K}^2))$ are two weak solutions for Eq.(\ref{eq1}) with $\hat{X}^1_{0}=\hat{X}^2_0$. If $\hat{\mP}\{(\hat{X}^1_t,\hat{K}^1_t)=(\hat{X}^2_t,\hat{K}^2_t), t\in[0,T]\}=1$, we say that the pathwise uniqueness holds for Eq.(\ref{eq1}).
\ed

In the following, we give some conditions to assure the existence and pathwise uniqueness of weak solutions for Eq.(\ref{eq1}). Assume:

\begin{enumerate}[($\bf{H}_{1.1}$)]
	\item The function $b$ is continuous in $(x,\mu)$, and $b,\sigma$ satisfy for $(x,\mu)\in\mR^{d}\times{\cP_{2}(\mR^d)}$
	\ce
	{|{b(x,\mu)}|}^2+\|\sigma(x,\mu)\|^2\leq{L_1(1+|{x}|^2+\mu(|\cdot|^2))},
	\de
	where $L_1>0$ is a constant.
\end{enumerate}

\begin{enumerate}[($\bf{H}_{1.2}$)]
	\item The functions $b,\sigma$ satisfy for $(x_1,\mu_1), (x_2,\mu_2)\in\mR^{d}\times{\cP_{2}(\mR^d)}$
	\ce
	&&2\langle{x_1-x_2,b(x_1,\mu_1)-b(x_2,\mu_2)}\rangle\leq {L_2\(|x_1-x_2|^2+\mW_2^2(\mu_1,\mu_2)\)},\\
	&&\parallel{\sigma(x_1,\mu_1)-\sigma(x_2,\mu_2)}\parallel^2\leq {L_2\(|x_1-x_2|^2+\mW_2^2(\mu_1,\mu_2)\)},
	\de
	where $L_2>0$ is a constant.
\end{enumerate}

Next, we give a key lemma. Set 
\ce
&&\sC^{\sL_{\xi}}_{0,T}:=\left\{\mu_{\cdot}\in C([0,T], \cP_2(\mR^d)), \mu_{0}=\sL_{\xi}\right\},\\
&&\hat{\rho}(\mu_{\cdot}, \nu_{\cdot}):=\sup\limits_{t\in[0,T]}\mW_2(\mu_t, \nu_t), \quad \mu_{\cdot}, \nu_{\cdot}\in \sC^{\sL_{\xi}}_{0,T},
\de
and the space $(\sC^{\sL_{\xi}}_{0,T}, \hat{\rho})$ is a complete metric space. For $\mu_{\cdot}\in \sC^{\sL_{\xi}}_{0,T}$, consider the following auxiliary multivalued SDE on $\mR^d$:
\be\left\{\begin{array}{l}
\dif X^{\mu_{\cdot}}_t\in \ -A(X^{\mu_{\cdot}}_t)\dif t+ \ b(X^{\mu_{\cdot}}_t,\mu_t)\dif t+\sigma(X^{\mu_{\cdot}}_t,\mu_t)\dif W_t,\\
X^{\mu_{\cdot}}_0=\xi. 
\end{array}
\label{meq1}
\right.
\ee

\bl\label{auxilemm}
Suppose that ${\rm Int}(\cD(A))\neq\emptyset$ and $b, \s$ satisfy $(\bf{H}_{1.1})$-$(\bf{H}_{1.2})$. Then for any $\sF_{0}$-measurable random variable $\xi$ with $\mE|\xi|^{2}<\infty$, Eq.(\ref{meq1}) has a unique strong solution $(X_{\cdot}^{\mu_{\cdot}},K_{\cdot}^{\mu_{\cdot}})$ and $\sL_{X^{\mu_{\cdot}}_{\cdot}}\in\sC^{\sL_{\xi}}_{0,T}$.
\el
\begin{proof}
First of all, by the similar deduction to that in \cite[Theorem 2.8]{rwz}, we obtain that Eq.(\ref{meq1}) has a unique strong solution $(X_{\cdot}^{\mu_{\cdot}},K_{\cdot}^{\mu_{\cdot}})$. Then, we prove that $\sL_{X^{\mu_{\cdot}}_{\cdot}}\in\sC^{\sL_{\xi}}_{0,T}$.

Take any $a\in {\rm Int}(\cD(A))$. By It\^{o}'s formula, Lemma \ref{inteineq} and $(\bf{H}_{1.1})$, we have for $0<t\leq T$
\be
\left| X_{t}^{\mu_{\cdot}}-a \right|^2&=&\left| \xi-a\right|^2-2\int_0^t\left< X_{s}^{\mu_{\cdot}}-a,\dif K_{s}^{\mu_{\cdot}} \right>  +2\int_0^t\< X_{s}^{\mu_{\cdot}}-a,b\left( X_{s}^{\mu_{\cdot}},\mu _s \right) \>\dif s\no\\
&&+2\int_0^t\left< X_{s}^{\mu_{\cdot}}-a,\sigma \left( X_{s}^{\mu_{\cdot}},\mu _s \right)\dif W_s \right>+\int_0^t\left\| \sigma \left( X_{s}^{\mu_{\cdot}},\mu _s \right) \right\| ^2\dif s\no\\
&\leq&\left| \xi-a\right|^2-2\gamma_1|K^{\mu_{\cdot}}|_0^t+2\gamma_2\int_0^t|X_{s}^{\mu_{\cdot}}-a|\dif s+2\gamma_3t+\int_0^t|X_{s}^{\mu_{\cdot}}-a|^2\dif s\no\\
&&+\int_0^t|b\left( X_{s}^{\mu_{\cdot}},\mu _s \right)|^2\dif s+\int_0^t\left\| \sigma \left( X_{s}^{\mu_{\cdot}},\mu _s \right) \right\| ^2\dif s\no\\
&&+2\int_0^t\left< X_{s}^{\mu_{\cdot}}-a,\sigma \left( X_{s}^{\mu_{\cdot}},\mu _s \right)\dif W_s \right>\no\\
&\leq&\left| \xi-a\right|^2-2\gamma_1|K^{\mu_{\cdot}}|_0^t+2\gamma_2\int_0^t(1+|X_{s}^{\mu_{\cdot}}-a|^2)\dif s+2\gamma_3t+\int_0^t|X_{s}^{\mu_{\cdot}}-a|^2\dif s\no\\
&&+L_1\int_0^t\(1+2|a|^2+2|X_{s}^{\mu_{\cdot}}-a|^2+\mu _s(|\cdot|^2)\)\dif s+2\int_0^t\left< X_{s}^{\mu_{\cdot}}-a,\sigma \left( X_{s}^{\mu_{\cdot}},\mu _s \right)\dif W_s \right>\no\\
&\leq&\left| \xi-a\right|^2-2\gamma_1|K^{\mu_{\cdot}}|_0^t+\(2\gamma_2+2\gamma_3+L_1(1+2|a|^2+\sup\limits_{s\in[0,T]}\mu _s(|\cdot|^2))\)T\no\\
&&+(2\gamma_2+1+2L_1)\int_0^t|X_{s}^{\mu_{\cdot}}-a|^2\dif s+2\int_0^t\left< X_{s}^{\mu_{\cdot}}-a,\sigma \left( X_{s}^{\mu_{\cdot}},\mu _s \right)\dif W_s \right>.\label{momeesti}
\ee
By the BDG inequality and the H\"older inequality, it holds that
\ce
&&\mathbb{E}\left( \underset{s\in \left[ 0,t \right]}{\sup}\left| X_{s}^{\mu_{\cdot}}-a \right|^2 \right)+2\gamma_1\mE|K^{\mu_{\cdot}}|_0^t\\
&\leqslant&\mathbb{E}\left| \xi-a \right|^2+\(2\gamma_2+2\gamma_3+L_1(1+2|a|^2+\sup\limits_{s\in[0,T]}\mu _s(|\cdot|^2))\)T\\
&&+(2\gamma_2+1+2L_1)\mE\int_0^t|X_{s}^{\mu_{\cdot}}-a|^2\dif s\\
&&+12\mE\(\int_0^t|X_{s}^{\mu_{\cdot}}-a|^2\|\sigma \left( X_{s}^{\mu_{\cdot}},\mu _s \right)\|^2\dif s\)^{1/2}\\
&\leqslant&\mathbb{E}\left| \xi-a \right|^2+\(2\gamma_2+2\gamma_3+L_1(1+2|a|^2+\sup\limits_{s\in[0,T]}\mu _s(|\cdot|^2))\)T\\
&&+(2\gamma_2+1+2L_1)\mE\int_0^t|X_{s}^{\mu_{\cdot}}-a|^2\dif s\\
&&+\frac{1}{2}\mathbb{E}\left( \underset{s\in \left[ 0,t \right]}{\sup}\left| X_{s}^{\mu_{\cdot}}-a \right|^2 \right)+C\mE\int_0^t\left\| \sigma \left( X_{s}^{\mu_{\cdot}},\mu _s \right) \right\| ^2\dif s,
\de
and furthermore
\ce
&&\mathbb{E}\left( \underset{s\in \left[ 0,t \right]}{\sup}\left| X_{s}^{\mu_{\cdot}}-a \right|^2 \right)+4\gamma_1\mE|K^{\mu_{\cdot}}|_0^t\\
&\leqslant&2\mathbb{E}\left| \xi-a \right|^2+2\(2\gamma_2+2\gamma_3+CL_1(1+2|a|^2+\sup\limits_{s\in[0,T]}\mu _s(|\cdot|^2))\)T\\
&&+2(2\gamma_2+1+2CL_1)\int_0^t\mathbb{E}\left( \underset{r\in \left[ 0,s\right]}{\sup}|X_{r}^{\mu_{\cdot}}-a|^2\right)\dif s.
\de
Then by the Gronwall inequality, we know that 
\ce
\mathbb{E}\left( \underset{s\in \left[ 0,t \right]}{\sup}\left| X_{s}^{\mu_{\cdot}}-a \right|^2 \right)\leq\[2\mathbb{E}\left| \xi-a \right|^2+C_{\gamma_2,\gamma_3,L_1}T\]e^{2(2\gamma_2+1+2CL_1)t},
\de
where $C_{\gamma_2,\gamma_3,L_1}:=2\(2\gamma_2+2\gamma_3+CL_1(1+2|a|^2+\sup\limits_{s\in[0,T]}\mu _s(|\cdot|^2))\)$, which yields that
\ce
\mathbb{E}\left( \underset{s\in \left[ 0,t \right]}{\sup}\left| X_{s}^{\mu_{\cdot}} \right|^2 \right)&\leq&2\mathbb{E}\left( \underset{s\in \left[ 0,t \right]}{\sup}\left| X_{s}^{\mu_{\cdot}}-a \right|^2 \right)+2|a|^2\no\\
&\leq&2\[2\mathbb{E}\left| \xi-a \right|^2+C_{\gamma_2,\gamma_3,L_1}T\]e^{2(2\gamma_2+1+2CL_1)t}+2|a|^2. 
\de
Thus, it holds that $\sL_{X_{t}^{\mu_{\cdot}}}\in\cP_2(\mR^d)$.

Next, we estimate $\mW_2(\sL_{X_{s}^{\mu_{\cdot}}}, \sL_{X_{t}^{\mu_{\cdot}}})$ for $s,t\in[0,T]$. From the definition of the metric $\mW_2$, it follows that
\be
\mW^2_2(\sL_{X_{s}^{\mu_{\cdot}}}, \sL_{X_{t}^{\mu_{\cdot}}})\leq \mE|X^{\mu_{\cdot}}_s-X^{\mu_{\cdot}}_t|^2.
\label{measesti}
\ee
So, it is sufficient to compute $\mE|X^{\mu_{\cdot}}_s-X^{\mu_{\cdot}}_t|^2$. 

Note that $X^{\mu_{\cdot}}_s, X^{\mu_{\cdot}}_t$ satisfy the following equations
\ce
&&X^{\mu_{\cdot}}_s=\xi-K^{\mu_{\cdot}}_s+ \int_0^s b(X^{\mu_{\cdot}}_r,\mu_r)\dif r+\int_0^s\sigma(X^{\mu_{\cdot}}_r,\mu_r)\dif W_r,\\
&&X^{\mu_{\cdot}}_t=\xi-K^{\mu_{\cdot}}_t+ \int_0^t b(X^{\mu_{\cdot}}_r,\mu_r)\dif r+\int_0^t\sigma(X^{\mu_{\cdot}}_r,\mu_r)\dif W_r.
\de
Assume $s\leq t$ and take a stopping time sequence $\{\tau_N\}$ given by $\tau_N:=\inf\{r\geq s, |X^{\mu_{\cdot}}_r|\geq N\}$. Thus, by the It\^o formula, it holds that
\ce
|X^{\mu_{\cdot}}_{t\land\tau_N}-X^{\mu_{\cdot}}_s|^2&=&-2\int_s^{t\land\tau_N}<X^{\mu_{\cdot}}_{r}-X^{\mu_{\cdot}}_s, \dif K^{\mu_{\cdot}}_r>+2\int_s^{t\land\tau_N}<X^{\mu_{\cdot}}_{r}-X^{\mu_{\cdot}}_s, b(X^{\mu_{\cdot}}_r,\mu_r)>\dif r\\
&&+2\int_s^{t\land\tau_N}<X^{\mu_{\cdot}}_{r}-X^{\mu_{\cdot}}_s, \sigma(X^{\mu_{\cdot}}_r,\mu_r)\dif W_r>+\int_s^{t\land\tau_N}\|\sigma(X^{\mu_{\cdot}}_r,\mu_r)\|^2\dif r\\
&\leq&\int_s^{t\land\tau_N}|X^{\mu_{\cdot}}_{r}-X^{\mu_{\cdot}}_s|^2\dif r+\int_s^{t\land\tau_N}|b(X^{\mu_{\cdot}}_r,\mu_r)|^2\dif r\\
&&+2\int_s^{t\land\tau_N}<X^{\mu_{\cdot}}_{r}-X^{\mu_{\cdot}}_s, \sigma(X^{\mu_{\cdot}}_r,\mu_r)\dif W_r>+\int_s^{t\land\tau_N}\|\sigma(X^{\mu_{\cdot}}_r,\mu_r)\|^2\dif r.
\de
Taking the expectation on two sides, by $(\bf{H}_{1.1})$ we obtain that
\ce
\mE|X^{\mu_{\cdot}}_{t\land\tau_N}-X^{\mu_{\cdot}}_s|^2&\leq&\mE\int_s^{t\land\tau_N}|X^{\mu_{\cdot}}_{r}-X^{\mu_{\cdot}}_s|^2\dif r+L_1\mE\int_s^{t\land\tau_N}(1+|X^{\mu_{\cdot}}_r|^2+\mu_r(|\cdot|^2))\dif r\\
&\leq&(1+2L_1)\mE\int_s^{t\land\tau_N}|X^{\mu_{\cdot}}_{r}-X^{\mu_{\cdot}}_s|^2\dif r+L_1\mE\int_s^{t\land\tau_N}(1+2|X^{\mu_{\cdot}}_s|^2+\mu_r(|\cdot|^2))\dif r\\
&\leq&(1+2L_1)\mE\int_s^{t}|X^{\mu_{\cdot}}_{r\land\tau_N}-X^{\mu_{\cdot}}_s|^2\dif r\\
&&+L_1\mE\(1+2|X^{\mu_{\cdot}}_s|^2+\sup\limits_{r\in[0,T]}\mu_r(|\cdot|^2)\)(t-s).
\de
Using the Gronwall inequality (\cite[Lemma 19.2, Page 172]{hzy}) and the mean value theorem, one can get that
\ce
\mE|X^{\mu_{\cdot}}_{t\land\tau_N}-X^{\mu_{\cdot}}_s|^2&\leq& L_1\mE\(1+2|X^{\mu_{\cdot}}_s|^2+\sup\limits_{r\in[0,T]}\mu_r(|\cdot|^2)\)\frac{e^{(1+2L_1)(t-s)}-1}{1+2L_1}\\
&\leq& L_1\mE\(1+2|X^{\mu_{\cdot}}_s|^2+\sup\limits_{r\in[0,T]}\mu_r(|\cdot|^2)\)(t-s)e^{(1+2L_1)t},
\de 
 and furthermore as $N\to\infty$
\be
\mE|X^{\mu_{\cdot}}_s-X^{\mu_{\cdot}}_t|^2\leq L_1\(1+2\mE|X^{\mu_{\cdot}}_s|^2+\sup\limits_{r\in[0,T]}\mu_r(|\cdot|^2)\)(t-s)e^{(1+2L_1)t}.
\label{diffesti}
\ee

Inserting (\ref{diffesti}) in (\ref{measesti}), we have that
\ce
\mW_2(\sL_{X_{s}^{\mu_{\cdot}}}, \sL_{X_{t}^{\mu_{\cdot}}})\leq L^{1/2}_1\(1+2\mE|X^{\mu_{\cdot}}_s|^2+\sup\limits_{r\in[0,T]}\mu_r(|\cdot|^2)\)^{1/2}(t-s)^{1/2}e^{(1+2L_1)t/2},
\de
which yields that 
$$
\lim\limits_{s\to t}\mW_2(\sL_{X_{s}^{\mu_{\cdot}}}, \sL_{X_{t}^{\mu_{\cdot}}})=0.
$$
By the same deduction to above, it holds that for $s\geq t$,
$$
\lim\limits_{s\to t}\mW_2(\sL_{X_{s}^{\mu_{\cdot}}}, \sL_{X_{t}^{\mu_{\cdot}}})=0.
$$
The proof is complete.
\end{proof}

Now, it is the position to state and prove the existence and uniqueness of strong solutions for Eq.(\ref{eq1}).

\bt\label{EU}
Assume that ${\rm Int}(\cD(A))\neq\emptyset$ and the coefficients $b, \sigma$ satisfy $(\bf{H}_{1.1})$-$(\bf{H}_{1.2})$. Then for any $\sF_{0}$-measurable random variable $\xi$ with $\mE|\xi|^{2}<\infty$, there exists a unique strong solution $(X,K)$ for Eq.(\ref{eq1}), i.e. for any $t\in[0,T]$
$$
X_{t}=\xi-K_{t}+\int_{0}^{t}b(X_{s},\sL_{X_{s}})\dif s +\int_{0}^{t}\sigma(X_{s},\sL_{X_{s}})\dif W_{s},\quad \sL_{X_{t}}=\mP\circ X_{t}^{-1}\in\cP_{2}(\mR^{d}).
$$
\et
\begin{proof}
{\bf Step 1.} Define the mapping ${\bf \Psi}: \sC^{\sL_{\xi}}_{0, t_0}\to \sC^{\sL_{\xi}}_{0, t_0}$ by ${\bf \Psi}(\mu_{\cdot})=\sL_{X^{\mu_{\cdot}}_{\cdot}}$ for any $\mu_{\cdot}\in\sC^{\sL_{\xi}}_{0, t_0}$, where $0<t_0\leq T$ is given later. Thus, by Lemma \ref{auxilemm}, we know that  the mapping ${\bf \Psi}$ is well defined. Then, we prove that it is contractive. 

First of all, by the definition of the metric on $\sC^{\sL_{\xi}}_{0, t_0}$, it holds that for $\mu^1_{\cdot}, \mu^2_{\cdot}\in\sC^{\sL_{\xi}}_{0, t_0}$, 
\be
\hat{\rho}\({\bf \Psi}(\mu^1_{\cdot}), {\bf \Psi}(\mu^2_{\cdot})\)&=&\sup\limits_{t\in[0,t_0]}\mW_2\({\bf \Psi}(\mu^1_{\cdot})_t, {\bf \Psi}(\mu^2_{\cdot})_t\)=\sup\limits_{t\in[0,t_0]}\mW_2\(\sL_{X^{\mu^1_{\cdot}}_t}, \sL_{X^{\mu^2_{\cdot}}_t})\no\\
&\leq&\sup\limits_{t\in[0,t_0]}\(\mE|X^{\mu^1_{\cdot}}_t-X^{\mu^2_{\cdot}}_t|^2\)^{1/2}\leq\(\mE\sup\limits_{t\in[0,t_0]}|X^{\mu^1_{\cdot}}_t-X^{\mu^2_{\cdot}}_t|^2\)^{1/2}.
\label{mappmetr1}
\ee

Next, we estimate $\mE\sup\limits_{t\in[0,t_0]}|X^{\mu^1_{\cdot}}_t-X^{\mu^2_{\cdot}}_t|^2$. By the It\^o formula, Lemma \ref{equi} and $(\bf{H}_{1.2})$, we have for $0<t\leq t_0$
\ce
\left| X^{\mu^1_{\cdot}}_t-X^{\mu^2_{\cdot}}_t\right|^2&=&-2\int_0^t\left< X^{\mu^1_{\cdot}}_s-X^{\mu^2_{\cdot}}_s,\dif \left( K_{s}^{\mu^1_{\cdot}}-K_{s}^{\mu^2_{\cdot}} \right) \right> \\
&&+2\int_0^t{\left< X^{\mu^1_{\cdot}}_s-X^{\mu^2_{\cdot}}_s,b\left( X^{\mu^1_{\cdot}}_s,\mu _{s}^{1} \right) -b\left( X^{\mu^2_{\cdot}}_s,\mu _{s}^{2} \right) \right>}\dif s\\
&&+2\int_0^t{\left< X^{\mu^1_{\cdot}}_s-X^{\mu^2_{\cdot}}_s,\(\sigma \left( X^{\mu^1_{\cdot}}_s,\mu _{s}^{1} \right) -\sigma \left( X^{\mu^2_{\cdot}}_s,\mu _{s}^{2} \right)\)\dif W_s \right>}\\
&&+\int_0^t{\left\| \sigma \left( X^{\mu^1_{\cdot}}_s,\mu _{s}^{1}  \right) -\sigma \left( X^{\mu^2_{\cdot}}_s,\mu _{s}^{2}\right) \right\|}^2\dif s\\
&\leq&2L_2\int_0^t\(\left| X^{\mu^1_{\cdot}}_s-X^{\mu^2_{\cdot}}_s \right|^2+\mW_2^2\left( \mu _{s}^{1},\mu _{s}^{2} \right)\)\dif s\\
&&+2\int_0^t{\left< X^{\mu^1_{\cdot}}_s-X^{\mu^2_{\cdot}}_s,\(\sigma \left( X^{\mu^1_{\cdot}}_s,\mu _{s}^{1}  \right) -\sigma \left( X^{\mu^2_{\cdot}}_s,\mu _{s}^{2}\right)\)\dif W_s \right>}
\de
So, it follows from the BDG inequality and the H\"older inequality that
\ce
\mathbb{E}\left( \underset{s\in \left[ 0,t \right]}{\sup}\left| X^{\mu^1_{\cdot}}_s-X^{\mu^2_{\cdot}}_s \right|^2 \right)&\leq&2L_2\mE\int_0^t\(\left| X^{\mu^1_{\cdot}}_s-X^{\mu^2_{\cdot}}_s \right|^2+\mW_2^2\left( \mu _{s}^{1},\mu _{s}^{2} \right)\)\dif s\no\\
&&+12\mE\left(\int_0^t| X^{\mu^1_{\cdot}}_s-X^{\mu^2_{\cdot}}_s|^2\left\|\sigma \left( X^{\mu^1_{\cdot}}_s,\mu _{s}^{1}  \right) -\sigma \left( X^{\mu^2_{\cdot}}_s,\mu _{s}^{2}\right)\right\|^2\dif s\right)^{1/2}\no\\ 
&\leq&2L_2\mE\int_0^t\(\left| X^{\mu^1_{\cdot}}_s-X^{\mu^2_{\cdot}}_s \right|^2+\mW_2^2\left( \mu _{s}^{1},\mu _{s}^{2} \right)\)\dif s\no\\
&&+C\int_0^t\mE\left\| \sigma \left( X^{\mu^1_{\cdot}}_s,\mu _{s}^{1}  \right) -\sigma \left( X^{\mu^2_{\cdot}}_s,\mu _{s}^{2}\right)\right\|^2\dif s\no\\
&&+\frac{1}{2}\mathbb{E}\left( \underset{s\in \left[ 0,t \right]}{\sup}\left| X^{\mu^1_{\cdot}}_s-X^{\mu^2_{\cdot}}_s \right|^2 \right),
\de
and furthermore
\be
\mathbb{E}\left( \underset{s\in \left[ 0,t \right]}{\sup}\left| X^{\mu^1_{\cdot}}_s-X^{\mu^2_{\cdot}}_s  \right|^2 \right)&\leq&4L_2\mE\int_0^t\(\left| X^{\mu^1_{\cdot}}_s-X^{\mu^2_{\cdot}}_s  \right|^2+\mW_2^2\left( \mu _{s}^{1},\mu _{s}^{2} \right)\)\dif s\no\\
&&+C\int_0^t\mE\left\| \sigma \left( X^{\mu^1_{\cdot}}_s,\mu _{s}^{1}  \right) -\sigma \left( X^{\mu^2_{\cdot}}_s,\mu _{s}^{2}\right)\right\|^2\dif s\no\\
&\overset{(\bf{H}_{1.2})}{\leq}&C\mE\int_0^t\(\left| X^{\mu^1_{\cdot}}_s-X^{\mu^2_{\cdot}}_s \right|^2+\mW_2^2\left( \mu _{s}^{1},\mu _{s}^{2} \right)\)\dif s\label{sequesti}\\
&\leq&Ct_0\sup\limits_{s\in[0,t_0]}\mW_2^2\left( \mu _{s}^{1},\mu _{s}^{2} \right)+C\int_0^t\mE\(\sup\limits_{r\in[0,s]}\left| X^{\mu^1_{\cdot}}_r-X^{\mu^2_{\cdot}}_r\right|^2\)\dif s.\no
\ee
By the Gronwall inequality, we derive that
\ce
\mathbb{E}\left( \underset{s\in \left[ 0,t_0 \right]}{\sup}\left| X^{\mu^1_{\cdot}}_s-X^{\mu^2_{\cdot}}_s \right|^2 \right) \leqslant Ct_0e^{Ct_0}\underset{s\in \left[ 0,t_0\right]}{\sup}\mW_2^2\left( \mu _{s}^{1},\mu _{s}^{2} \right) .
\de
Taking $t_{0}>0$ with $Ct_0e^{Ct_0}<\frac{1}{4}$, one can obtain that 
\be
\left( \mathbb{E}\left( \underset{s\in \left[ 0,t_0 \right]}\sup\left| X^{\mu^1_{\cdot}}_s-X^{\mu^2_{\cdot}}_s \right|^2 \right)\right)^{1/2} &\leqslant \frac{1}{2}\underset{s\in \left[ 0,t_0 \right]}{\sup}\mW_2\left( \mu _{s}^{1},\mu _{s}^{2} \right)=\frac{1}{2}\hat{\rho}(\mu^1_{\cdot}, \mu^2_{\cdot}). 
\label{mappmetr2}
\ee
Combining (\ref{mappmetr1}) with (\ref{mappmetr2}), we have 
\ce
\hat{\rho}\({\bf \Psi}(\mu^1_{\cdot}), {\bf \Psi}(\mu^2_{\cdot})\)\leq \frac{1}{2}\hat{\rho}(\mu^1_{\cdot}, \mu^2_{\cdot}). 
\de

{\bf Step 2.} We prove that Eq.(\ref{eq1}) has a unique strong solution.

By the conclusion in {\bf Step 1.} and the fixed point theorem, it holds that there exists a unique fixed point $\mu^*_{\cdot}\in\sC^{\sL_{\xi}}_{0, t_0}$ such that 
$$
{\bf \Psi}(\mu^*_{\cdot})=\sL_{X^{\mu^*_{\cdot}}_{\cdot}}=\mu^*_{\cdot}.
$$
Thus, $(X^{\mu^*_{\cdot}}_{\cdot}, K^{\mu^*_{\cdot}}_{\cdot})$ is a weak solution for Eq.(\ref{eq1}) on $[0, t_0]$. If $t_0\geq T$, the proof of the existence for weak solutions to Eq.(\ref{eq1}) is over; if $t_0<T$, by the same deduction as above, we obtain the existence for weak solutions to Eq.(\ref{eq1}) on $[t_0, T]$.

Next, we prove the pathwise uniqueness for Eq.(\ref{eq1}). Assume that $(\Omega, \mathscr{F}, \{\mathscr{F}_t\}_{t\in[0,T]}, \mP; W,\\ X_{\cdot}, K_{\cdot} ) $ and $(\Omega, \mathscr{F}, \{\mathscr{F}_t\}_{t\in[0,T]}, \mP; W, \tilde{X}_{\cdot}, \tilde{K}_{\cdot})$ are two weak solutions of Eq.(\ref{eq1}) with $X_0=\tilde{X}_0=\xi$, i.e.
\ce
&&X_{t}=\xi-K_{t}+\int_{0}^{t}b(X_{s},\sL_{X_{s}})\dif s +\int_{0}^{t}\sigma(X_{s},\sL_{X_{s}})\dif W_{s},\\
&&\tilde{X}_{t}=\xi-\tilde{K}_{t}+\int_{0}^{t}b(\tilde{X}_{s},\sL_{\tilde{X}_{s}})\dif s +\int_{0}^{t}\sigma(\tilde{X}_{s},\sL_{\tilde{X}_{s}})\dif W_{s}.
\de
So, by the similar deduction to that in (\ref{sequesti}), we obtain that for any $t\in[0,T]$
\ce
\mathbb{E}\left( \underset{s\in \left[ 0,t \right]}{\sup}\left| X_{s}-\tilde{X}_{s} \right|^2\right)&\leq& C\mE\int_0^t\(\left| X_{s}-\tilde{X}_{s} \right|^2+\mW_2^2\left( \sL_{X_{s}}, \sL_{\tilde{X}_{s}}\right)\)\dif s\\
&\leq&C\mE\int_0^t\(\left| X_{s}-\tilde{X}_{s} \right|^2+\mE\left| X_{s}-\tilde{X}_{s} \right|^2\)\dif s\\
&\leq&C\int_0^t\mathbb{E}\left( \underset{r\in \left[ 0,s\right]}{\sup}\left| X_{r}-\tilde{X}_{r} \right|^2\right)\dif s
\de
It follows from the Gronwall inequality that
$$
\mathbb{E}\left( \underset{s\in \left[ 0,t \right]}{\sup}\left| X_{s}^{}-\tilde{X}_{s} \right|^2 \right) =0,
$$
which yields that
$$
X_{t}=\tilde{X}_t, \quad t\in[0,T] \,\, a.s.\mathbb{P}.
$$

Finally, note that for any $t\in[0,T]$
\ce
K_t&=&\xi-X_t+\int_0^t{b\left( X_{s},\sL_{X_{s}} \right) \dif s+\int_0^t{\sigma \left( X_{s},\sL_{X_{s}}\right) \dif W_s}}\\
&=&\xi-\tilde{X_t}+\int_0^t{b\left(\tilde{X}_{s},\sL_{\tilde{X}_{s}} \right) \dif s+\int_0^t{\sigma \left( \tilde{X}_{s},\sL_{\tilde{X}_{s}}\right) \dif W_s}}\\
&=&\tilde{K_t}.           
\de
Thus, the fact that $K_t$ is continuous in $t$ assures that $K_{t}=\tilde{K}_t, t\in[0,T] \,\, a.s.\mathbb{P}$. The proof is complete.
\end{proof}

By Theorem \ref{EU}, we immediately have the following corollary.

\bc
Assume that ${\rm Int}(\cD(A))\neq\emptyset$ and the coefficients $b, \sigma$ satisfy $(\bf{H}_{1.1})$-$(\bf{H}_{1.2})$. Then for any $\sF_{0}$-measurable random variable $\xi$ with $\mE|\xi|^{2}<\infty$, the strong solution $(X,K)$ for Eq.(\ref{eq1}) has the following moment property: for $0\leq t\leq T$
\be
\mathbb{E}\left( \underset{s\in \left[ 0,t \right]}{\sup}\left| X_{s}\right|^2\right)\leq 2(2\mathbb{E}\left| \xi-a \right|^2+Ct)e^{Ct}+2|a|^2, \quad \forall a\in {\rm Int}(\cD(A)). 
\label{itoieq2}  
\ee
\ec

\section{The generalized It\^o formula for multivalued McKean-Vlasov SDEs}\label{geneito}

In the section, we will extend the classical It\^{o} formula for SDEs  to multivalued McKean-Vlasov SDEs. First of all, we strengthen the condition ($\bf{H}_{1.2}$) to the following assumption:
\begin{enumerate}[($\bf{H}_{1.3}$)]
	\item The functions $b,\sigma$ satisfy for $(x_1,\mu_1), (x_2,\mu_2)\in\mR^{d}\times{\cP_{2}(\mR^d)}$
	\ce
	&&\left| b(x_1,\mu_1)-b(x_2,\mu_2)\right|^2 +\parallel{\sigma(x_1,\mu_1)-\sigma(x_2,\mu_2)}\parallel^2\\
	&\leq&{L_3\(|x_1-x_2|^2+\mW_2^2(\mu_1,\mu_2)\)},
	\de
	where $L_3>0$ is a constant.
\end{enumerate}

\bt\label{if}
Suppose that the function $\Phi:\mR^d\times\cP_{2}(\mR^d)\mapsto\mR$ belongs to $C_{b}^{2,2}(\mR^d\times\cP_{2}(\mR^d))$. Then, under ${\rm Int}(\cD(A))\neq\emptyset$ and $(\bf{H}_{1.1})$, $(\bf{H}_{1.3})$, for any $\sF_{0}$-measurable random variable $\xi$ with $\mE|\xi|^{2}<\infty$, the following It\^{o} formula holds:\,\,$\forall 0\leq s <t $,
\be
&&\Phi(X_t,\sL_{X_{t}})-\Phi(X_s,\sL_{X_{s}})\no\\
&=&-\int_s^t\left< \partial _x\Phi\left( X_u,\sL_{X_{u}} \right),\dif K_u \right> +\int_s^t(b^i\partial_{x_i}\Phi)(X_u,\sL_{X_{u}})\dif u\no\\
&&+\int_s^t(\sigma^{ij}\partial_{x_i} \Phi)(X_u,\sL_{X_{u}})\dif W_u^j+\frac{1}{2}\int_s^t\((\sigma\sigma^*)^{ij}\partial_{x_ix_j}^2\Phi\)(X_u,\sL_{X_{u}})\dif u\no\\
&&+\int_s^t\int_{\mR^d}b^i(y,\sL_{X_{u}})(\partial_\mu \Phi)_i(X_u,\sL_{X_{u}})(y)\sL_{X_{u}}(\dif y)\dif u\no\\
&&+\frac{1}{2}\int_s^t\int_{\mR^d}(\sigma\sigma^*)^{ij}(y,\sL_{X_{u}})\partial_{y_i}(\partial_\mu \Phi)_j(X_u,\sL_{X_{u}})(y)\sL_{X_{u}}(\dif y)\dif u\no\\
&&-\int_s^t\mE\left<\left( \partial _{\mu}\Phi \right) \left( y,\sL_{X_{u}} \right) \left( X_{u} \right),\dif K_u \right>|_{y=X_u}.
\label{ito1}
\ee
\et
\begin{proof}
{\bf Step 1.} Suppose that $b, \sigma$ are bounded. Fix $x\in\mR^d$ and define $f(\mu):=\Phi(x,\mu) $. Now, we prove the It\^o formula for $f(\sL_{X_{t}})$.

For any positive integer $N$, set
\be
x^1, x^2,\cdots, x^N\in\mR^d, \quad f^N(x^1, x^2,\cdots, x^N):=f\(\frac{1}{N}\sum_{l=1}^N\delta_{x^l}\),
\label{hkde}
\ee
and $f^N(x^1, x^2,\cdots, x^N)$ is a function on $\mR^{d\times N}$. Moreover, by \cite[Proposition 3.1, Page 15]{ccd}, it holds that $f^N$ is $C^2$ on $\mR^{d\times N}$ and 
\be
\partial_{x^i}f^N(x^1, x^2,\cdots, x^N)&=&\frac{1}{N}\partial_{\mu}f\(\frac{1}{N}\sum_{l=1}^N\delta_{x^l}\)(x^i),\no\\
\partial^2_{x^ix^j}f^N(x^1, x^2,\cdots, x^N)&=&\frac{1}{N}\partial_{x}\partial_{\mu}f\(\frac{1}{N}\sum_{l=1}^N\delta_{x^l}\)(x^i)I_{i=j}+\frac{1}{N^2}\partial^2_{\mu}f\(\frac{1}{N}\sum_{l=1}^N\delta_{x^l}\)(x^i, x^j),\no\\
\qquad\qquad i,j=1,2,\cdots, N.
\label{deri}
\ee
Besides, we take $N$ independent copies $X^l_t, l=1,2,\cdots, N$ of $X_t$. That is, consider these following equations
\ce\left\{\begin{array}{l}
\dif X_{t}^{l}\in-A(X_{t}^{l})\dif t+b(X_{t}^{l},\sL_{X_{t}^{l}})\dif t+\sigma(X_{t}^{l},\sL_{X_{t}^{l}})\dif W_{t}^{l},\\
X_{0}^{l}=\xi,  \qquad l=1,2,\cdots, N,
\end{array}
\right.
\de
where $W_{\cdot}^l, l=1,2,\cdots, N$ are mutually independent and identical distribution copies of $W_{\cdot}$. By Theorem \ref{EU}, we know that for $l=1,2,\cdots, N$, there exists $(X_{\cdot}^l,K_{\cdot}^l)\in \sA$ such that 
$$
X_{t}^{l}-X_{s}^{l}=-(K_{t}^l-K_{s}^l)+\int_s^tb(X_{r}^{l},\sL_{X_{r}^{l}})\dif r+\int_s^t\sigma(X_{r}^{l},\sL_{X_{r}^{l}})\dif W_{r}^l, \quad 0\leq s<t.
$$
Then applying It\^{o}'s formula to $f^N\left( X_{t}^{1},\cdots ,X_{t}^{N} \right) $ and taking the expectation on both sides, we obtain that for $0\leq t<t+v$
\ce
&&\mathbb{E}f^N\left( X_{t+v}^{1},\cdots ,X_{t+v}^{N} \right)\\
&=&\mathbb{E}f^N\left( X_{t}^{1},\cdots ,X_{t}^{N} \right)+\sum_{i=1}^N{\int_t^{t+v}{\mathbb{E}\partial _{x^i}f^N\left( X_{s}^{1},\cdots ,X_{s}^{N} \right)}b\left( X_{s}^{i},\sL_{X_{s}^{i}}\right)\dif s}
\\
&&+\frac{1}{2}\sum_{i=1}^N{\int_t^{t+v}{\mathbb{E}\partial _{x^ix^i}^{2}f^N\left( X_{s}^{1}, \cdots, X_{s}^{N} \right)}\sigma \sigma ^*\left( X_{s}^{i},\sL_{X_{s}^{i}}\right) \dif s}
\\
&&-\sum_{i=1}^N\mathbb{E}{\int_t^{t+v}{\left< \partial _{x^i}f^N\left( X_{s}^{1},\cdots ,X_{s}^{N} \right),\dif K_{s}^{i} \right>}}
\\
&=&\mathbb{E}f^N\left( X_{t}^{1},\cdots ,X_{t}^{N} \right) +N\int_t^{t+v}{\mathbb{E}\partial _{x^1}f^N\left( X_{s}^{1},\cdots ,X_{s}^{N} \right)}b\left( X_{s}^{1},\sL_{X_{s}^{1}} \right) \dif s
\\
&&+\frac{N}{2}\int_t^{t+v}{\mathbb{E}\partial _{x^1x^1}^{2}f^N\left( X_{s}^{1},\cdots ,X_{s}^{N} \right)}\sigma \sigma ^*\left( X_{s}^{1},\sL_{X_{s}^{1}} \right) \dif s
\\
&&-N\mathbb{E}\int_t^{t+v}{ \left<\partial _{x^1}f^N\left( X_{s}^{1},\cdots ,X_{s}^{N} \right),\dif K_{s}^{1} \right>},
\de
where the property of the same distributions for $(X^{l}_t, K^{l}_t), l=1,2,\cdots, N$ is used in the second equality. Inserting (\ref{hkde}) and (\ref{deri}) in the above equality, we get that
\ce
\mathbb{E}f\left( \frac{1}{N}\sum_{l=1}^N{\delta _{X_{t+v}^{l}}} \right) &=&\mathbb{E}f\left( \frac{1}{N}\sum_{l=1}^N{\delta _{X_{t}^{l}}} \right) +\int_t^{t+v}{\mathbb{E}\partial _{\mu}f\left( \frac{1}{N}\sum_{l=1}^N{\delta _{X_{s}^{l}}} \right)}\left( X_{s}^{1} \right) b\left( X_{s}^{1},\sL_{X_{s}^{1}}  \right) \dif s\\
&&+\frac{1}{2}\int_t^{t+v}{\mathbb{E}\partial _y\partial _{\mu}f\left(\frac{1}{N}\sum_{l=1}^N{\delta _{X_{s}^{l}}}\right)}\left( X_{s}^{1} \right) \sigma \sigma ^*\left( X_{s}^{1},\sL_{X_{s}^{1}}  \right) \dif s\\
\\
&&+\frac{1}{2N}\int_t^{t+v}{\mathbb{E}\partial _{\mu}^{2}f\left( \frac{1}{N}\sum_{l=1}^N{\delta _{X_{s}^{l}}}\right)}\left( X_{s}^{1},X_{s}^{1} \right) \sigma \sigma ^*\left( X_{s}^{1},\sL_{X_{s}^{1}} \right) \dif s
\\
&&-\mathbb{E}\int_t^{t+v}{\left< \partial _{\mu}f\left( \frac{1}{N}\sum_{l=1}^N{\delta _{X_{s}^{l}}}\right) \left( X_{s}^{1} \right) ,\dif K_{s}^{1} \right>}.	
\de

Next, we take the limit on both sides of the above equality. Note that 
$$
\underset{N\rightarrow \infty}{\lim}\mathbb{E}\left[\mW_2^2\left( \frac{1}{N}\sum_{l=1}^N{\delta _{X_{t}^{l}}},\sL_{X_{t}}  \right) \right] =0,
$$
which comes from \cite[Section 5]{HK}. And as $N\to \infty$, by continuity and boundedness of $f, \partial_{\mu}f, \partial_{y}\partial_{\mu}f$, and boundedness of $\partial^{2}_{\mu}f, b, \sigma$, it follows from the dominated convergence theorem that 
\be
f\left( \sL_{X_{t+v}}  \right) &=&f\left( \sL_{X_{t}}  \right) +\int_t^{t+v}{\mathbb{E}\partial _{\mu}f\left( \sL_{X_{s}}  \right) \left( X_{s}^{1} \right)}b\left( X_{s}^{1},\sL_{X^1_{s}} \right) \dif s\no\\
&&+\frac{1}{2}\int_t^{t+v}{\mathbb{E}\partial _y\partial _{\mu}f\left( \sL_{X_{s}}  \right) \left( X_{s}^{1} \right)}\sigma \sigma ^*\left( X_{s}^{1},\sL_{X^1_{s}} \right) \dif s\no \\
&&-\mathbb{E}\int_t^{t+v}{\left< \partial _{\mu}f\left( \sL_{X_{s}}  \right) \left( X_{s}^{1} \right),\dif K_{s}^{1} \right>}\no \\
&=&f\left( \sL_{X_{t}}  \right) +\int_t^{t+v}{\int_{\mR^d}\partial _{\mu}f\left( \sL_{X_{s}}  \right) \left( y\right)}b\left( y,\sL_{X_{s}} \right) \sL_{X_{s}}(\dif y)\dif s\no\\
&&+\frac{1}{2}\int_t^{t+v}{\int_{\mR^d}\partial _y\partial _{\mu}f\left( \sL_{X_{s}}  \right) \left(y\right)}\sigma \sigma ^*\left( y,\sL_{X_{s}} \right)\sL_{X_{s}}(\dif y) \dif s\no \\
&&-\int_t^{t+v}{\mathbb{E}\left< \partial _{\mu}f\left( \sL_{X_{s}}  \right) \left( X_{s} \right),\dif K_{s}\right>}.
\label{feq2}
\ee

{\bf Step 2.} Assume that $(\bf{H}_{1.1})$ $(\bf{H}_{1.3})$ hold. Then we prove the It\^o formula for $f(\sL_{X_{t}})$.

Let $\phi_{n}:\mR^d \mapsto \mR^d$ be a smooth function satisfying 
\begin{equation*}
	\begin{cases}
		\phi_{n}(x)=x, & \left| x\right|\leq n,\\
		\phi_{n}(x)=0, & \left| x\right|> 2n,
	\end{cases}
\end{equation*}
\\
such that 
\be
\left|\phi_{n}(x) \right|\leq C, \quad \left| \left| \partial_x \phi_{n}(x)\right| \right| \leq C,
\label{auxiboun}
\ee
where the positive constant $C$ is independent of $n$. Define 
$$
b^{\left( n \right)}\left( x,\mu \right) :=b\left( \phi _n\left( x \right) ,\mu \right), \quad \sigma ^{\left( n \right)}\left( x,\mu \right) :=\sigma \left( \phi _n\left( x \right) ,\mu \right),
$$
and by simple calculation it holds that $b^{\left( n \right)}, \sigma ^{\left( n \right)} $ satisfy the assumption $(\bf{H}_{1.1})$-$(\bf{H}_{1.2})$. Therefore, the following equation
\be\left\{\begin{array}{l}
\dif X_{t}^{\left( n\right) }\in-A(X_{t}^{\left( n\right) })\dif t+b^{\left( n\right) }(X_{t}^{\left( n\right) },\sL_{X_{t}^{\left( n\right) }})\dif t+\sigma^{\left( n\right) }(X_{t}^{\left( n\right) },\sL_{X_{t}^{\left( n\right) }})\dif W_{t},\\
X_{0}^{\left( n\right)}=\xi,
\end{array}
\right.
\label{seq1}
\ee
has a unique solution $(X_{\cdot}^{\left( n\right) },K_{\cdot}^{\left( n\right) })\in \sA$. Besides, by $(\bf{H}_{1.1})$, we know that $b^{\left( n \right)}, \sigma ^{\left( n \right)} $ are bounded. Thus, by {\bf Step 1.}, it holds that for $0\leq t<t+v$,
\be
f\left( \sL_{X_{t+v}^{\left( n\right) }} \right) &=&f\left( \sL_{X_{t}^{\left( n\right) }} \right) +\int_t^{t+v}{\int_{\mR^d}\partial _{\mu}f\left( \sL_{X^{(n)}_{s}}  \right) \left( y\right)}b^{\left( n \right)}\left( y,\sL_{X^{(n)}_{s}} \right) \sL_{X^{(n)}_{s}}(\dif y)\dif s\no\\
&&+\frac{1}{2}\int_t^{t+v}{\int_{\mR^d}\partial _y\partial _{\mu}f\left( \sL_{X^{(n)}_{s}}  \right) \left( y\right)}\sigma^{\left( n \right)}(\sigma^{\left( n \right)})^*\left( y,\sL_{X^{(n)}_{s}} \right)\sL_{X^{(n)}_{s}}(\dif y) \dif s\no \\
&&-\int_t^{t+v}{\mathbb{E}\left< \partial _{\mu}f\left( \sL_{X^{(n)}_{s}}  \right) \left( X^{(n)}_{s} \right),\dif K^{(n)}_{s}\right>}.
\label{feq12}
\ee

Next, we observe the limit of $\sL_{X_{t}^{\left( n\right) }}$ as $n\rightarrow \infty$. Applying It\^{o}'s formula to $\left| X_{t}^{\left(n \right) } -X_{t}\right|^{2} $ and taking the expectation on both sides, we get that
\ce
\mathbb{E}\left| X_{t}^{\left( n \right)}-X_t \right|^2&=&-2\mathbb{E}\int_0^t{\left< X_{s}^{\left( n \right)}-X_s,\dif K_{s}^{\left( n \right)}-\dif K_s \right>}\\
&&+2\mathbb{E}\int_0^t{\left< X_{s}^{\left( n \right)}-X_s,\left( b_{}^{\left( n \right)}\left( X_{s}^{\left( n \right)},\sL_{X_{s}^{\left( n \right)}} \right) -b\left( X_{s},\sL_{X _{s}}\right) \right) \right>}\dif s\\
&&+\mathbb{E}\int_0^t{\left\| \sigma _{}^{\left( n \right)}\left( X_{s}^{\left( n \right)},\sL_{X_{s}^{\left( n \right)}}\right) -\sigma \left( X_{s}, \sL_{X _{s}}\right) \right\|}^2\dif s\\
&\leq&\mathbb{E}\int_0^t{\left| X_{s}^{\left( n \right)}-X_s \right|}^2\dif s+\mathbb{E}\int_0^t\left|b_{}^{\left( n \right)}\left( X_{s}^{\left( n \right)},\sL_{X_{s}^{\left( n \right)}} \right) -b\left( X_{s},\sL_{X _{s}}\right)\right|^2\dif s\\
&&+\mathbb{E}\int_0^t{\left\| \sigma _{}^{\left( n \right)}\left( X_{s}^{\left( n \right)},\sL_{X_{s}^{\left( n \right)}}\right) -\sigma \left( X_{s}, \sL_{X _{s}}\right) \right\|}^2\dif s\\
&\leqslant& \mathbb{E}\int_0^t{\left| X_{s}^{\left( n \right)}-X_s \right|}^2\dif s+L_3\mathbb{E}\int_0^t{\left( \left| \phi _n\left( X_{s}^{\left( n \right)} \right) -X_s \right|^2 +\mW_2^2\left( \sL_{X_{s}^{\left( n \right)}}, \sL_{X _{s}} \right) \right)}\dif s
\\
&\leqslant& C\mathbb{E}\int_0^t{\left| X_{s}^{\left( n \right)}-X_s \right|}^2\dif s+C\mathbb{E}\int_0^T\left| \phi _n\left( X_{s} \right) -X_s \right|^2 \dif s,
\de
where we use $(\bf{H}_{1.3})$, (\ref{auxiboun}) and the fact $\mW_2^2\left( \sL_{X_{s}^{\left( n \right)}}, \sL_{X _{s}} \right) \leqslant \mathbb{E}\left| X_{s}^{\left( n \right)}-X_s \right|^2$. Therefore, by the Gronwall inequality, it holds that 
$$
\mathbb{E}\left| X_{t}^{\left( n \right)}-X_t \right|^2\leqslant C\mathbb{E}\int_0^T \left| \phi _n\left( X_{s}\right) -X_s \right|^2 \dif s.
$$
Combining $\left| \phi _n\left( x \right) \right|\leqslant C$ for $x\in \mR^d$, with the estimate (\ref{itoieq2}), by the dominated convergence theorem one can have that
$$
\underset{n\rightarrow \infty}{\lim}\mathbb{E}\int_0^T{\left| \phi _n\left( X_{s}^{} \right) -X_s \right|^2}\dif s=0,
$$
and
$$
\underset{n\rightarrow \infty}{\lim}\mathbb{E}\left| X_{t}^{\left( n \right)}-X_t \right|^2=0.
$$
Moreover, we obtain that 
\be
\underset{n\rightarrow \infty}{\lim}\mW_2^2\left( \sL_{X_{t}^{\left( n \right)}}, \sL_{X _{t}} \right)\leqslant \underset{n\rightarrow \infty}{\lim}\mathbb{E}\left| X_{t}^{\left( n \right)}-X_t \right|^2=0.
\label{statdistesti}
\ee

Next, note that $(X_{\cdot}^{(n)}, K_{\cdot}^{(n)})$, $(X_{\cdot}, K_{\cdot})$ are the unique strong solution of Eq.(\ref{seq1}) and (\ref{eq1}), respectively, i.e.
\ce
&&K_{t}^{(n)}=\xi-X_{t}^{(n)}+\int_0^t{b^{\left( n\right) }\left( X_{s}^{(n)},\sL_{ X_{s}^{(n)}} \right)\dif s+\int_0^t{\sigma^{\left( n\right) } \left( X_{s}^{(n)},\sL_{ X_{s}^{(n)}}  \right)\dif W_s}},\\
&&K_t=\xi-X_t+\int_0^t{b\left( X_s,\sL_{X_s} \right)\dif s+\int_0^t{\sigma \left( X_s,\sL_{X_s}  \right)\dif W_s}}.
\de
Thus by (\ref{statdistesti}) and $(\bf{H}_{1.3})$, we get 
$$
\underset{n\rightarrow \infty}{\lim}\mathbb{E}\left| K_{t}^{\left( n \right)}-K_t \right|^2=0.
$$
Taking the limit on two sides of (\ref{feq12}), by the dominated convergence theorem, one can conclude (\ref{feq2}).

{\bf Step 3.} We prove the It\^{o} formula (\ref{ito1}). 

By the classical It\^{o}'s formula and (\ref{feq2}), it holds that
\ce
&&\Phi(X_{t+v},\sL_{X _{t+v}})-\Phi(X_t,\sL_{X _{t}}) \\
&=&-\int_t^{t+v}{\left< \partial _x\Phi\left( X_s,\sL_{X _{s}} \right),\dif K_s \right>}+\int_t^{t+v}{\left<  \right. \partial _x}\Phi\left( X_s,\sL_{X _{s}} \right) ,b\left( X_s,\sL_{X _{s}}\right) \left.  \right> \dif s\\
&&+\int_t^{t+v}{\left< \partial _x\Phi\left( X_s,\sL_{X _{s}} \right) ,\sigma \left( X_s,\sL_{X _{s}} \right) \dif W_s \right>}+\frac{1}{2}\int_t^{t+v}{{\rm tr}\left( \sigma \sigma ^*\left( X_s,\sL_{X _{s}} \right) \partial _{x}^{2}\Phi\left( X_s,\sL_{X _{s}} \right) \right)}\dif s\\
&&+\int_t^{t+v}{\int_{\mathbb{R}^d}^{}{\left< b\left( y,\sL_{X _{s}} \right) ,\partial _{\mu}\Phi\left( X_s,\sL_{X _{s}} \right) \left( y \right) \right> \sL_{X _{s}}\left( \dif y \right)}}\dif s\\
&&+\frac{1}{2}\int_t^{t+v}{\int_{\mathbb{R}^d}{{\rm tr}\left( \sigma \sigma ^*\left( y, \sL_{X _{s}} \right) \partial _y\partial _{\mu}\Phi\left( X_s,\sL_{X _{s}} \right) \left( y \right) \right)}}\sL_{X _{s}}\left( \dif y \right) \dif s\\
&&-\int_t^{t+v}{\mE\left<\partial _{\mu}\Phi\left( y,\sL_{X _{s}} \right) \left( X _{s} \right),\dif K_s \right>|_{y=X_s}}.
\de
The proof is complete.
\end{proof}

\br
Although we prove the It\^o formula for any $\Phi\in C_{b}^{2,2}(\mR^d\times\cP_{2}(\mR^d))$, it also holds for $\Phi\in \sC(\mR^d\times\cP_{2}(\mR^d))$. Please refer to \cite[Proposition A.8]{WDL} for details.
\er

\section{The stability of strong solutions}\label{stab}

In the section, we require that $\xi=x_0$ is non-random and  study the stability of strong solutions for Eq.(\ref{eq1}) by the generalized It\^o formula.   

\subsection{The asymptotic stability of the second moment for the strong solution}\label{eub}

In the subsection, we consider the asymptotic stability of the second moment for the strong solution of Eq.$(\ref{eq1})$. Assume:

\begin{enumerate}[($\bf{H}_{2.1}$)]
	\item There exists a function $F:\mR^d\times\cP_{2}(\mR^d)\mapsto\mR$ satisfying\\
	$(i)$ $F\in \sC(\mR^d\times\cP_{2}(\mR^d))$,\\
	$(ii)$
	$$
	\int_{\mR^d} \( \cL_\mu F(x,\mu)+\alpha F(x,\mu) \) \mu(\dif x)\leq M_1,
	$$
	where $\cL_\mu$ is defined by
	\ce
	\(\cL_\mu F\)(x,\mu)&:=&\(b^i\partial_{x_i} F\)(x,\mu)+\frac{1}{2}\((\sigma\sigma^*)^{ij}\partial_{x_ix_j}^2 F\)(x,\mu)\\
	&&+\int_{\mR^d}b^i(y,\mu)(\partial_\mu F)_i(x,\mu)(y)\mu(\dif y)\\
	&&+\int_{\mR^d}\frac{1}{2}(\sigma\sigma^*)^{ij}(y,\mu)\partial_{y_i}(\partial_\mu F)_j(x,\mu)(y)\mu(\dif y),
	\de
	and $\alpha>0, M_1\geq0$ are constants; \\
	$(iii)$ 
	$$
	a_1\int_{\mR^d}|x|^2\mu(\dif x)-M_2\leq \int_{\mR^d} F(x,\mu)\mu(\dif x),
	$$
	where $a_1>0, M_2\geq0$ are constants.
\end{enumerate}

In the following, we prepare an important lemma.
\bl\label{averlimi}
Assume that ${\rm Int}(\cD(A))\neq\emptyset$ and $(\bf{H}_{1.1})$ $(\bf{H}_{1.3})$ and $(\bf{H}_{2.1})$ hold. If the strong solution $(X_{\cdot},K_{\cdot})$ and the Lyapunov function $F$ satisfy for any $t>0$
\be
\left< \partial _xF\left( X_t,\sL_{X _{t}} \right), \dif K_t \right> +\mE\left<{\left( \partial _{\mu}F\right) \left( x,\sL_{X _{t}} \right) \left( X_t\right)},\dif K_t \right>|_{x=X_t}\geq 0,
\label{strolyap}
\ee
it holds that
\ce
\limsup\limits_{t\rightarrow\infty}\frac{1}{t}\int_0^t\mE|X _{s}|^2\dif s\leq \frac{\a M_2+M_1}{\a a_1}.
\de
Moreover, if $M_1=M_2=0$, 
\be
\int_0^{\infty}\mE|X _{s}|^2\dif s<\infty.
\label{inte}
\ee
\el
\begin{proof}
By the It\^{o} formula (\ref{ito1}), it holds that
\ce
&&F(X_t,\sL_{X _{t}})-F(x_0,\d_{x_0})\\
&=&\int_0^t\[(b^i\partial_{x_i}F)(X_s,\sL_{X _{s}})+\frac{1}{2}\((\sigma\sigma^*)^{ij}\partial_{x_ix_j}^2F\)(X_s,\sL_{X _{s}})\no\\
&&\qquad+\int_{\mR^d}b^i(y,\sL_{X _{s}})(\partial_\mu F)_i(X_s,\sL_{X _{s}})(y)\sL_{X _{s}}(\dif y)\\
&&\qquad+\frac{1}{2}\int_{\mR^d}(\sigma\sigma^*)^{ij}(y,\sL_{X _{s}})\partial_{y_i}(\partial_\mu F)_j(X_s,\sL_{X _{s}})(y)\sL_{X _{s}}(\dif y)\]\dif s\no\\
&&-\int_0^t\(\left< \partial _xF\left( X_s,\sL_{X _{s}} \right), \dif K_s \right> +\mE\left< \left(\partial _{\mu}F \right) \left( x,\sL_{X _{s}} \right) \left( X_s \right),\dif K_s \right>|_{x=X_s}\)\\
&&+\int_0^t(\sigma^{ij}\partial_{x_i}F)(X_s,\sL_{X _{s}})\dif W_s^j\no\\
&\overset{(\ref{strolyap})}{\leq}&\int_0^t\cL_\mu F(X_s,\sL_{X _{s}})\dif s+\int_0^t(\sigma^{ij}\partial_{x_i}F)(X_s,\sL_{X _{s}})\dif W_s^j.
\de
Localizing and taking the expectation on both sides of the above equality, by the Fatou lemma and $(ii), (iii)$ of $(\bf{H}_{2.1})$ we obtain that
\ce
\mE F(X_t,\sL_{X _{t}})-F(x_0,\d_{x_0})&\leq& \mE\int_0^t\cL_\mu F(X_s,\sL_{X _{s}})\dif s=\int_0^t\int_{\mR^d}\cL_\mu F(x,\sL_{X _{s}})\sL_{X _{s}}(\dif x)\dif s\\
&\leq& -\a\int_0^t\int_{\mR^d} F(x,\sL_{X _{s}})\sL_{X _{s}}(\dif x)\dif s+M_1 t\\
&\leq&-\a \int_0^t\(a_1\int_{\mR^d}|x|^2\sL_{X _{s}}(\dif x)-M_2\)\dif s+M_1 t\\
&=&-\a a_1\int_0^t\mE|X _{s}|^2\dif s+(\a M_2+M_1)t,
\de
which concludes that
\ce
\mE F(X_t,\sL_{X _{t}})\leq F(x_0,\d_{x_0})-\a a_1\int_0^t\mE|X _{s}|^2\dif s+(\a M_2+M_1)t.
\de
Finally, based on $(iii)$ in $(\bf{H}_{2.1})$, one can get that
\be
a_1\mE|X_t|^2-M_2\leq F(x_0,\d_{x_0})-\a a_1\int_0^t\mE|X _{s}|^2\dif s+(\a M_2+M_1)t,
\label{iiiu}
\ee
which yields that
\be
\frac{1}{t}\int_0^t\mE|X _{s}|^2\dif s\leq \frac{F(x_0,\d_{x_0})+M_2}{\a a_1 t}+\frac{\a M_2+M_1}{\a a_1},
\label{aver}
\ee
and furthermore
$$
\limsup\limits_{t\rightarrow\infty}\frac{1}{t}\int_0^t\mE|X _{s}|^2\dif s\leq \frac{\a M_2+M_1}{\a a_1}.
$$
The proof is complete.
\end{proof}

Here we state and prove the main result in the subsection.
\bt\label{asymstab}
Suppose that ${\rm Int}(\cD(A))\neq\emptyset$ and $(\bf{H}_{1.1})$ $(\bf{H}_{1.3})$ and $(\bf{H}_{2.1})$ hold with $M_1=M_2=0$. If the strong solution $(X_{\cdot},K_{\cdot})$ and the Lyapunov function $F$ satisfy for any $t>0$
\ce
\left< \partial _xF\left( X_t,\sL_{X _{t}} \right), \dif K_t \right> +\mE\left<{\left( \partial _{\mu}F\right) \left( x,\sL_{X _{t}} \right) \left( X_t\right)},\dif K_t \right>|_{x=X_t}\geq 0,
\de
it holds that 
\ce
\lim\limits_{t\rightarrow\infty}\mE|X_t|^2=0. 
\de
\et
\begin{proof}
By (\ref{inte}) of Lemma \ref{sp}, in order to obtain $\lim\limits_{t\rightarrow\infty}\mE|X_t|^2=0$, we only need to prove that $\mE|X_t|^2$ is uniformly continuous in $t$. Note that for any $s,t$
\ce
|\mE|X_t|^2-\mE|X_s|^2|\leq \mE|X_t-X_s|^2+2(\mE|X_t-X_s|^2)^{1/2}(\mE|X_s|^2)^{1/2}.
\de
Therefore it is sufficient to show that $\mE|X_t-X_s|^2$ uniformly converges to $0$ as $s\rightarrow t$. For $\mE|X_t-X_s|^2$, by the same deduction to that of (\ref{diffesti}), we have that
\ce
\mE|X_t-X_s|^2\leq L_1\(1+3\mE|X_s|^2\)(t-s)e^{(1+2L_1)t},
\de
which together with (\ref{itoieq2}) and the continuity of $X_t$ in $t$ implies that $\mE|X_t-X_s|^2$ uniformly converges to $0$ as $s\rightarrow t$. The proof is complete.
\end{proof}

\br
Here we remind that one can prove the above theorem without (\ref{inte}). In fact, when $M_1=M_2=0$, (\ref{iiiu}) becomes 
$$
a_1\mE|X_t|^2\leq F(x_0,\d_{x_0})-\a a_1\int_0^t\mE|X _{s}|^2\dif s.
$$
By the comparison theorem for integral equations, it holds that
\ce
\mE|X_t|^2\leq\frac{F(x_0,\d_{x_0})}{a_1}e^{-\a t},
\de
which yields that
\ce
\lim\limits_{t\rightarrow\infty}\mE|X_t|^2=0. 
\de
\er

\subsection{The almost surely asymptotic stability for the strong solution}\label{asas}

In the subsection, we study the almost surely asymptotic stability of the strong solution for Eq.$(\ref{eq1})$. We start with the concept of the almost surely asymptotic stability.

\bd\label{asymstab}
If for all $x_0\in\mR^d$, it holds that
\ce
\mP\left\{\lim_{t\rightarrow\infty}|X_t|=0\right\}=1,
\de
we say that $X_{\cdot}$ is almost surely asymptotically stable.
\ed

Next we introduce a function class. Let $\varSigma$ denote the family of functions $\gamma:\mR_+\mapsto\mR_+$, which are continuous, strictly increasing, and $\gamma(0)=0$. And $\varSigma_\infty$ means the family of functions $\gamma\in \varSigma$ with $\gamma(x)\rightarrow\infty$ as $x\rightarrow\infty$. Then we present some assumptions.

\begin{enumerate}[($\bf{H}^{\prime}_{1.1}$)]
\item The function $b$ is continuous in $(x,\mu)$ and satisfies for $(x,\mu)\in\mR^{d}\times{\cP_2(\mR^d)}$
\ce
{|{b(x,\mu)}|}^2\leq{L'_1(1+|{x}|^2+\mu(|\cdot|^2))},
\de
where $L'_1>0$ is a constant, and $\sigma$ is bounded.
\end{enumerate}

\begin{enumerate}[($\bf{H}_{2.2}$)]
\item There exists a function $F:\mR^d\times\cP_2(\mR^d)\mapsto\mR$ satisfying\\
$(i)$ $F\in C_{b,+}^{2,2;1}(\mR^d\times\cP_2(\mR^d))$,\\
$(ii)$
$\cL_\mu F(x,\mu)+\alpha F(x,\mu)\leq 0$,\\
$(iii)$
$\gamma_1(|x|)\leq F(x,\mu)\leq \gamma_2(|x|)$,
where $\gamma_1, \gamma_2\in\varSigma_\infty$.
\end{enumerate}

\bt\label{sp}
Assume that ${\rm Int}(\cD(A))\neq\emptyset$ and $(\bf{H}^{\prime}_{1.1})$ $(\bf{H}_{1.3})$ and $(\bf{H}_{2.2})$ hold. If the strong solution $(X_{\cdot},K_{\cdot})$ and the Lyapunov function $F$ satisfy for any $t>0$
\ce
\left< \partial _xF\left( X_t,\sL_{X _{t}} \right), \dif K_t \right> +\mE\left<{\left( \partial _{\mu}F\right) \left( x,\sL_{X _{t}} \right) \left( X_t\right)},\dif K_t \right>|_{x=X_t}\geq 0,
\de
$X_{\cdot}$ is almost surely asymptotically stable, i.e.
\ce
\mP\left\{\lim_{t\rightarrow\infty}|X_t|=0\right\}=1.
\de
\et
\begin{proof}
First of all, since under $(\bf{H}^{\prime}_{1.1})$ $(\bf{H}_{1.3})$ Eq.(\ref{eq1}) has a unique strong solution $(X_{\cdot},K_{\cdot})$ with the initial value $(x_0,0)$, the distribution family $\{\sL_{X _{t}}\}_{t\geq 0}$ of $(X_t)_{t\geq 0}$ is known. Thus, we rewrite Eq.(\ref{eq1}) as
\be
X_t=x_0-K_t+\int_0^t\tilde{b}(u, X_u)\dif u+\int_0^t\tilde{\sigma}(u, X_u)\dif W_u, \quad t\geq 0,
\label{Eq2}
\ee
where $\tilde{b}(u, X_u):=b(X_u,\sL_{X _{u}}), \tilde{\sigma}(u, X_u):=\sigma(X_u,\sL_{X _{u}})$. That is, Eq.(\ref{Eq2}) is a nonhomogeneous multivalued SDE. Set $\tau_n:=\inf\{t\geq0,|X_t|>n\}.$
Now applying It\^{o}'s formula to $|X_{s\wedge\tau_n}-x_0|^2$ for $s\geq0$, we get that
\ce
|X_{s\wedge\tau_n}-x_0|^2&=&-2\int_0^{s\wedge\tau_n}<X_u-x_0, \dif K_u>+2\int_0^{s\wedge\tau_n}{<X_u-x_0, \tilde{b}(u, X_u)>}\dif u\\
&&+2\int_0^{s\wedge\tau_n}<X_u-x_0, \tilde{\sigma}(u, X_u)\dif W_u>+\int_0^{s\wedge\tau_n}\|\tilde{\sigma}(u, X_u))\|^2\dif u\\
&\leq&2\int_0^{s\wedge\tau_n}{|X_u-x_0| | \tilde{b}(u, X_u)|}\dif u+\int_0^{s\wedge\tau_n}\|\tilde{\sigma}(u, X_u)\|^2\dif u\\
&&+2\int_0^{s\wedge\tau_n}<X_u-x_0, \tilde{\sigma}(u, X_u)\dif W_u>\\
&\leq&\int_0^{s\wedge\tau_n}|X_u-x_0|^2\dif u+\int_0^{s\wedge\tau_n}|\tilde{b}(u, X_u)|^2\dif u+\int_0^{s\wedge\tau_n}\|\tilde{\sigma}(u, X_u)\|^2\dif u\\
&&+2\int_0^{s\wedge\tau_n}<X_u-x_0, \tilde{\sigma}(u, X_u)\dif W_u>.\\
\de
Then by $(\bf{H}^{\prime}_{1.1})$ and the Burkholder-Davis-Gundy inequality, it holds that
\ce
\mE\left(\sup_{0\leq s\leq t}|X_{s\wedge\tau_n}-x_0|^2\right)&\leq&\mE\int_0^{t\wedge\tau_n}|X_u-x_0|^2\dif u+C\mE\int_0^{t\wedge\tau_n}(1+| X_u|^2+n^2)\dif u\\
&&+C\mE\left(\int_0^{t\wedge\tau_n}|X_u-x_0|^2\|\tilde{\sigma}(u, X_u)\|^2\dif u\right)^{1/2}\\
&\leq&\mE\int_0^{t\wedge\tau_n}|X_u-x_0|^2\dif u+C\mE\int_0^{t\wedge\tau_n}(1+|X_u|^2+n^2)\dif u\\
&&+\frac{1}{4}\mE\left(\sup_{0\leq s\leq t}|X_{s\wedge\tau_n}-x_0|^2\right)+C\mE\left(\int_0^{t\wedge\tau_n}\|\tilde{\sigma}(u, X_u)\|^2\dif u\right),
\de
which yields
\ce
\mE\left(\sup_{0\leq s\leq t}|X_{s\wedge\tau_n}-x_0|^2\right)\leq C\mE\int_0^{t\wedge\tau_n}|X_u-x_0|^2du+C\mE\int_0^{t\wedge\tau_n}(1+| X_u|^2+n^2)\dif u.
\de
So, based on the boundedness of $X_t$ we conclude that
\ce
\mE\left(\sup_{0\leq s\leq t}|X_{s\wedge\tau_n}-x_0|^2\right)\leq C\mE({t\wedge\tau_n})\leq Ct,
\de
where $C>0$ depends on $L'_1$, $x_0$ and $n$. Then from the Chebyshev inequality, it follows that for any $\lambda>0$, 
\be
\mP\left\{\sup_{0\leq s\leq t}|X_{s\wedge\tau_n}-x_0|>\lambda\right\}\leq\frac{Ct}{\lambda^2}.\label{estiprob}
\ee

Next, we follow up the line in \cite[Theorem 5.2]{DQ2} and apply (\ref{estiprob}) to obtain that
\ce
\mP\left\{\lim_{t\rightarrow\infty}|X_t|=0\right\}=1.
\de
The proof is complete.
\end{proof}

\section{An example}\label{app}

In this section, we present an example to explain our results. 

\bx
Let $\cO:=\left\{x \in \mathbb{R}^d ; l(x) \geqslant 0\right\}$ be a convex closed domain with $\partial \cO=\left\{x \in \mathbb{R}^d ; l(x)=0\right\}$, where $l$ belongs to $C_b^2\left(\mathbb{R}^d\right)$ with $\sum\limits_{i=1}^d(\p_{x_i}l(x))^2=1$ and $\<x+\frac{1}{4}\int_{\mR^d}y\mu(\dif y),\\ \p_x l(x))\>\geq 0$ for $x\in\partial \cO$ and $\mu\in\cP_2(\partial \cO)$.

Next, we consider the following equation
\be\left\{\begin{array}{l}
\dif X_t\in \ -\partial I_{\cO}(X_t)\dif t+\(-X_t-\frac{1}{4}\int_{\mR^d}y\sL_{X_t}(\dif y)\)\dif t+\dif W_t,\\
X_0=x_0\in \bar\cO,
\end{array}
\label{xeq}
\right.
\ee
where 
$$
A=\partial I_{\cO}, \quad b(x,\mu)=-x-\frac{1}{4}\int_{\mR^d}y\mu(\dif y), \quad \sigma(x,\mu)=I_d.
$$
Moreover, by {\bf Example} \ref{consetex}, it holds that
\ce\left\{\begin{array}{l}
\dif X_t=\p_x l(X_t)\dif |K|^t_0+\(-X_t-\frac{1}{4}\int_{\mR^d}y\sL_{X_t}(\dif y)\)\dif t+\dif W_t,\\
X_0=x_0\in \bar\cO, \quad  {\rm supp}(\dif |K|^t_0)\subset \{u\geq 0: X_u\in \p \cO\}.
\end{array}
\label{exeq}
\right.
\de
So, one can justify that for $(x,\mu)\in\mR^d\times{\cP_2(\mR^d)}$
\ce
|b(x,\mu)|^2&=&\left|x+\frac{1}{4}\int_{\mR^d}y\mu(dy)\right|^2\leq2\(|x|^2+\frac{1}{16}\left|\int_{\mR^d}y\mu(dy)\right|^2\)\\
&\leq&2\(|x|^2+\frac{1}{16}\int_{\mR^d}|y|^2\mu(dy)\)\leq\frac{17}{8}\(1+|x|^2+\mu(|\cdot|^2)\),
\de
and for $(x_1,\mu_1), (x_2,\mu_2)\in\mR^{d}\times{\cP_2(\mR^d)}$, 
\ce
|b(x_1,\mu_1)-b(x_2,\mu_2)|^2&=&\left|x_1+\frac{1}{4}\int_{\mR^d}y\mu_1(dy)-x_2-\frac{1}{4}\int_{\mR^d}y\mu_2(dy)\right|^2\\
&\leq&2|x_1-x_2|^2+2\frac{1}{16}\left|\int_{\mR^d}y\mu_1(dy)-\int_{\mR^d}y\mu_2(dy)\right|^2\\
&\leq&\frac{17}{8}\(|x_1-x_2|^2+\mW_2^2(\mu_1,\mu_2)\).
\de
That is, $b, \sigma$ satisfy $(\bf{H}_{1.1})$ and $(\bf{H}_{1.3})$. By Theorem \ref{EU}, we know that Eq.(\ref{xeq}) has a unique strong solution $(X,K)$.

In the following, if we take the Lyapunov function 
$$
F(x,\mu)=\left|x+\frac{1}{4}\int_{\mR^d}y\mu(dy)\right|^2,
$$
$F$ satisfies $(\bf{H}_{2.1})$. Indeed, it is easily seen that $F$ satisfies $(i)$. And by simple calculation, it holds that
\ce
&&\p_x F(x,\mu)=2(x+\frac{1}{4}\int_{\mR^d}y\mu(dy)),\quad \p_\mu F(x,\mu)(z)=\frac{1}{2}(x+\frac{1}{4}\int_{\mR^d}y\mu(dy)),\\
&&\p_{xx} F(x,\mu)=2, \qquad\qquad \p_z\p_\mu F(x,\mu)=0.
\de
Based on this, the following computations are reasonable:
\ce
&&(\cL_\mu F)(x,\mu)\\
&=&-2\left|x+\frac{1}{4}\int_{\mR^d}y\mu(dy)\right|^2+d-\frac{1}{2}\int_{\mR^d}\(z_i+\frac{1}{4}\int_{\mR^d}y_i\mu(dy)\)\(x_i+\frac{1}{4}\int_{\mR^d}y_i\mu(dy)\)\mu(dz),
\de
and 
\ce
&&\int_{\mR^d}\(\cL_\mu F(x,\mu)+2 F(x,\mu)\)\mu(dx)\\
&=&-2\int_{\mR^d}\left|x+\frac{1}{4}\int_{\mR^d}y\mu(dy)\right|^2\mu(dx)+d-\sum_{i=1}^d\frac{1}{2}\left(\int_{\mR^d}\(x_i+\frac{1}{4}\int_{\mR^d}y_i\mu(dy)\)\mu(dx)\right)^2\\
&&+2\int_{\mR^d} \left|x+\frac{1}{4}\int_{\mR^d}y\mu(dy)\right|^2\mu(dx)\\
&\leq&d.
\de
That is, $F(x,\mu)$ satisfies $(ii)$. Besides, note that
\ce
F(x,\mu)&=&\sum_{i=1}^d\left(x_i+\frac{1}{4}\int_{\mR^d}y_i\mu(dy)\right)^2\geq \sum_{i=1}^d\left(x^2_i+\frac{x_i}{2}\int_{\mR^d}y_i\mu(dy)+\frac{1}{16}\left(\int_{\mR^d}y_i\mu(dy)\right)^2\right)\\
&\geq&\sum_{i=1}^d\left(\frac{3x^2_i}{4}-\frac{3}{16}\left(\int_{\mR^d}y_i\mu(dy)\right)^2\right)\geq \sum_{i=1}^d\left(\frac{3x^2_i}{4}-\frac{3}{16}\int_{\mR^d}y^2_i\mu(dy)\right).
\de
Thus, we have
\ce
\int_{\mR^d}F(x,\mu)\mu(dx)\geq \frac{9}{16}\int_{\mR^d}|x|^2\mu(dx),
\de
which implies that $F(x,\mu)$ satisfies $(iii)$. 

Now, we verify that (\ref{strolyap}) holds. In fact, note that $\<x+\frac{1}{4}\int_{\mR^d}y\mu(\dif y),\p_x l(x))\>\geq 0$ for $x\in\partial \cO$ and $\mu\in\cP_2(\partial \cO)$. Thus, it holds that $\<\p_x F(x,\mu),\p_x l(x))\>+\<\p_\mu F(x,\mu)(z),\p_x l(x))\>\geq 0$ for $x\in\partial \cO$ and $\mu\in\cP_2(\partial \cO)$. Therefore, (\ref{strolyap}) holds.

Finally, by the above deduction and Lemma \ref{averlimi}, we know that the strong solution of Eq.(\ref{xeq}) satisfies
$$
\limsup\limits_{t\rightarrow\infty}\frac{1}{t}\int_0^t\mE|X _{s}|^2\dif s\leq \frac{8}{9}d.
$$
\ex

\end{document}